\documentclass[a4paper,12pt,final]{amsart}
\usepackage{times,a4wide,mathrsfs,amssymb,amsmath,amsthm,enumerate,xypic,tikzsymbols,dsfont}

\newcommand{\C}{\mathbb{C}}

\newcommand{\LLL}{\mathbb{L}}
\newcommand{\QQ}{\mathbb{Q}}
\newcommand{\NN}{\mathbb{N}}
\newcommand{\PP}{\mathbb{P}}

\newcommand{\OO}{\mathcal O}

\newcommand{\XX}{\mathcal X}
\newcommand{\YY}{\mathcal Y}

\newcommand{\HH}{\mathcal H}
\newcommand{\NNN}{\mathcal N}
\newcommand{\VV}{\mathcal V}
\newcommand{\WW}{\mathcal W}

\newcommand{\Zz}{\mathcal Z}

\newcommand{\MM}{\mathcal M}

\newcommand{\FF}{\mathcal F}

\newcommand{\coker}{\hbox{Coker}}
\newcommand{\codim}{\hbox{codim}}

\newcommand{\wt}{\widetilde}
\newcommand{\rom}{\romannumeral}

\DeclareMathOperator{\Grif}{Grif}

\newcommand{\one}{\mathds{1}}

\newcommand{\hookdownarrow}{\mathrel{\rotatebox[origin=t]{90}{$\hookleftarrow$}}}

\newcommand\undermat[2]{
  \makebox[0pt][l]{$\smash{\underbrace{\phantom{
    \begin{matrix}#2\end{matrix}}}_{\text{$#1$}}}$}#2}

\DeclareMathOperator{\ima}{Im}

\DeclareMathOperator{\rank}{rank}
\DeclareMathOperator{\sym}{Sym}
\DeclareMathOperator{\Gr}{Gr}
\DeclareMathOperator{\Pf}{Pf}
\DeclareMathOperator{\Kuz}{Kuz}
\DeclareMathOperator{\Sing}{Sing}

\DeclareMathOperator{\Dec}{Dec}

\newtheorem{theorem}{Theorem}[section]

\newtheorem{lemma}[theorem]{Lemma}

\newtheorem{corollary}[theorem]{Corollary}
\newtheorem{proposition}[theorem]{Proposition}
\newtheorem{conjecture}[theorem]{Conjecture}
\newtheorem{convention}{Conventions}

\newtheorem{nonumbering}{Theorem}

\newtheorem{nonumberingc}{Corollary}

\theoremstyle{definition}
\newtheorem{remark}[theorem]{Remark}
\newtheorem{definition}[theorem]{Definition}

\newtheorem{example}[theorem]{Example}
\newtheorem{notation}[theorem]{Notation}

\newtheorem{nonumberingt}{Acknowledgments}

\begin{document}

\author[Robert Laterveer]
{Robert Laterveer}

\address{Institut de Recherche Math\'ematique Avanc\'ee,
CNRS -- Universit\'e 
de Strasbourg,\
7 Rue Ren\'e Des\-car\-tes, 67084 Strasbourg CEDEX,
FRANCE.}
\email{robert.laterveer@math.unistra.fr}

\title{Motives and the Pfaffian--Grassmannian equivalence}

\begin{abstract} We consider the Pfaffian--Grassmannian equivalence from the motivic point of view. The main result is that under certain numerical conditions, both sides of the equivalence are related on the level of Chow motives. The consequences include a verification of Orlov's conjecture for Borisov's Calabi--Yau threefolds, 
and verifications of Kimura's finite-dimensionality conjecture, Voevodsky's smash conjecture and the Hodge conjecture for certain linear sections of Grassmannians. We also obtain new examples of Fano varieties with infinite-dimensional Griffiths group.
 \end{abstract}


\thanks{Supported by ANR grant ANR-20-CE40-0023.}

\maketitle

\section{Introduction}

Given $V$ a complex vector space of dimension $n$, let
  \[   \Gr(2,V)\ \ \ \subset\ \PP(\wedge^2 V)   \]
  denote the Grassmannian of 2-dimensional linear subspaces of $V$ in its Pl\"ucker embedding.
  The projective dual to $ \Gr(2,V)$ is the {\em Pfaffian\/} of degenerate skew forms
  \[ \Pf:= \Bigl\{ \omega\in \PP(\wedge^2 V^\vee)\, \Big\vert\,  \hbox{rank} \, \omega < r_{max} \Bigr\}\ \ \ \subset\  \PP(\wedge^2 V^\vee)\ , \]
 where $r_{max}=n$ if $n$ is even, and $r_{max}=n-1$ if $n$ is odd. 
  
 Given a linear subspace $U\subset \wedge^2 V$ of codimension $k$, one can define varieties by intersecting on the Grassmannian side and on the Pfaffian side:
    \[  \begin{split}   X&=X_U:=  \Gr(2,V)\cap   \PP(U) \ \ \ \subset\ \PP(\wedge^2 V) \ , \\
                              Y&=Y_U:= \Pf\,\cap\, \PP(U^\perp)\ \ \ \subset\ \PP(\wedge^2 V^\vee) \ .
                              \end{split}  \]  
  For $U$ generic and $k$ small enough (the precise condition is that $k\le 6$ when $n$ is even, and $k\le 10$ when $n$ is odd, so that $Y$ avoids the singular locus of $\Pf$), the intersections $X$ and $Y$ are smooth and dimensionally transverse, of dimension
   \[  \begin{split}  \dim X&=2(n-2)-k\ ,\\
                             \dim Y&=\begin{cases} k-2 &\hbox{if\ $n$\ even},\\
                                                                  k-4 &\hbox{if\ $n$\ odd}.\\    
                                                                  \end{cases} \\
                                                                  \end{split} \]
                                                                  
These varieties have been intensively studied, and particularly so in the case $n=k=7$ \cite{Rod}, \cite{Bor}, \cite{BC}, \cite{Kuz0}, \cite{Mar}. In this case, $X$ and $Y$ are Calabi--Yau threefolds
that are L-equivalent \cite{Bor}, \cite{Mar} and derived equivalent \cite{BC}, \cite{Kuz0}, \cite{ADS}, while for general $U$ they are {\em not\/} birational. The L-equivalence of these Calabi--Yau threefolds gave rise to the first example that the affine line is a zero-divisor in the Grothendieck ring of varieties \cite{Bor}.

For $(n,k)$ arbitrary, the varieties $X$ and $Y$ can still be related on the level of the Grothendieck ring of varieties, and hence on the level of cohomology (cf.  subsections \ref{ss:l} and \ref{ss:coh} below), but the relation on the level of derived categories becomes conjectural (this is because the Pfaffian is singular, and one needs to find a categorical resolution of singularities, cf. subsection \ref{ss:der} below). The main result of the present paper relates $X$ and $Y$ on the level of Chow motives:

\begin{nonumbering}[=Theorem \ref{main}] Let $X$ and $Y$ be as above smooth dimensionally transverse intersections.
Assume $k<6$ is odd, or $(n,k)=(7,7)$. Assume also that the transcendental cohomology $H^\ast_{tr}(X,\QQ)$ is non-zero. Then there is an isomorphism of Chow motives
   \[ t(X)\ \xrightarrow{\cong}\ t(Y)(-m)\ \ \ \hbox{in}\ \MM_{\rm rat}\ ,\]
   where $m={1\over 2}(\dim X -\dim Y)$. 
\end{nonumbering}

Here $t(X)$ is a certain motive with the property that $h(X)=t(X)\oplus \bigoplus\one(\ast)$. The above is a simplified version; the actual statement of Theorem \ref{main} applies to certain other values of $(n,k)$.
Theorem \ref{main} means that the Grassmannian--Pfaffian equivalence is now better understood on the level of Chow motives than on the level of derived categories (which is as it should be: morally speaking, motives are easier to handle than derived categories).

In proving Theorem \ref{main}, we rely on the ``spread'' argument crafted by Voisin \cite{V0}, \cite{V1}, \cite{Vo}. This argument (used here in the form of the Franchetta-type result Proposition \ref{Frtype}) consists in working with families of varieties and correspondences, and exploiting the fact that the total space of the family has a very simple structure. Actually, the
starting point of the present paper was the realization that the Grassmannian--Pfaffian equivalence (and more generally, much of the set-up of HPD) is ideally suited to Voisin's ``spread'' argument.

The particular case $(n,k)=(7,7)$ of Theorem \ref{main} has the following consequence:

\begin{nonumberingc}[=Corollary \ref{CY3}] Let $X$ and $Y$ be the Calabi-Yau threefolds arising from the Pfaffian--Grassmannian equivalence \cite{Bor}.
Then there is an isomorphism of Chow motives
  \[ h(X)\cong h(Y)\ \ \ \hbox{in}\ \MM_{\rm rat}\ .\]
  \end{nonumberingc}
  
  This is in agreement with Orlov's conjecture \cite{Or}, stating that derived equivalent varieties should have isomorphic Chow motive.

For other values of $(n,k)$, Theorem \ref{main} can be applied to the study of Chow groups $A^i(X):=CH^i(X)_{\QQ}$ (i.e. the groups of codimension $i$ algebraic cycles on $Y$ with $\QQ$-coefficients, modulo rational equivalence). Here is a sampling of some applications:

\begin{nonumberingc}[=Corollaries \ref{12}, \ref{findim}, \ref{smash} and \ref{hc}]
Let $X$ be a smooth dimensionally transverse intersection
   \[ X:= \Gr(2,n)\cap H_1\cap\cdots\cap H_k\ \ \ \subset\ \PP^{{n\choose 2}-1} \ ,\]
   where the $H_j$ are hyperplanes.
   
  \noindent
  (\rom1) Assume $k\le 2$. Then
     \[  A^j_{hom}(X)=0\ \ \ \forall\ j\ .\]   
 \noindent
 (\rom2)  
  Assume $n$ is even and $k=3$, or $n$ is odd and $k=5$.
   Then 
   \[  A^j_{AJ}(X)=0\ \ \ \forall\ j\ ,\]
   and in particular
   $X$ has Kimura finite-dimensional motive. (Here $A^j_{AJ}(X)\subset A^j(X)$ denotes the subgroup of Abel--Jacobi trivial cycles.)

\noindent
(\rom3) Assume $n$ is even and $k=4$, or $n$ is odd and $k=6$.
   Then 
   \[  A^j_{AJ}(X)=0\ \ \ \forall\ j\not={1\over 2} \dim X +1\ ,\]
   and in particular
   Voevodsky's smash conjecture is true for $X$.
   
\noindent
(\rom4) Assume $n$ is even and $k=5$, or $n$ is odd and $k=6$.
Then 
  \[ A^j_{hom}(X)=0\ \ \ \forall \ j> {1\over 2}(\dim X+3)\ ,\]
and in particular the Hodge conjecture is true for $X$. (Here $A^j_{hom}(X)\subset A^j(X)$ denotes the subgroup of homologically trivial cycles.)
   \end{nonumberingc}
   
   There is also an application to the {\em Griffiths group\/} (i.e. homologically trivial algebraic cycles modulo algebraic equivalence):
   
 \begin{nonumberingc}[=Corollary \ref{grif}] Let $X$ be a general complete intersection
   \[ X:= \Gr(2,10)\cap H_1\cap\cdots\cap H_5\ \ \ \subset\ \PP^{44} \ ,\]
   where the $H_j$ are hyperplanes.
 Then $X$ is a Fano elevenfold and the Griffiths group $\Grif^6(X)_{\QQ}$ is infinite-dimensional.  
  \end{nonumberingc}

   It would be interesting to try and extend the results of this paper to generalized Pfaffian varieties (cf. \cite[Conjecture 6]{Kuz0} and \cite[Section 5.3]{Thom} and \cite{RS} for the conjectural HPD statement, and \cite{BorL} for a relation on the level of Hodge numbers).

   
 \vskip0.6cm

\begin{convention} In this article, the word {\sl variety\/} will refer to a reduced irreducible scheme of finite type over $\C$. A {\sl subvariety\/} is a (possibly reducible) reduced subscheme which is equidimensional. 

{\bf All Chow groups will be with rational coefficients}: we denote by $A_j(Y):=CH_j(Y)_{\QQ} $ the Chow group of $j$-dimensional cycles on $Y$ with $\QQ$-coefficients; for $Y$ smooth of dimension $n$ the notations $A_j(Y)$ and $A^{n-j}(Y)$ are used interchangeably. 
The notations $A^j_{hom}(Y)$ and $A^j_{AJ}(X)$ will be used to indicate the subgroup of homologically trivial (resp. Abel--Jacobi trivial) cycles.

The contravariant category of Chow motives (i.e., pure motives with respect to rational equivalence as in \cite{Sc}, \cite{MNP}) will be denoted 
$\MM_{\rm rat}$.
\end{convention}

\section{Preliminaries}

 \subsection{Cayley's trick and motives}

\begin{theorem}[Jiang \cite{Ji}]\label{ji} Let $ E\to U$ be a vector bundle of rank $r\ge 2$ over a smooth projective variety $U$, and let $S:=s^{-1}(0)\subset U$ be the zero locus of a regular section $s\in H^0(U,E)$ such that $S$ is smooth of dimension $\dim U-\rank E$. Let $X:=w^{-1}(0)\subset \PP(E)$ be the zero locus of the regular section $w\in H^0(\PP(E),\OO_{\PP(E)}(1))$ that corresponds to $s$ under the natural isomorphism $H^0(U,E)\cong H^0(\PP(E),\OO_{\PP(E)}(1))$, and assume $X$ is smooth. There is an isomorphism of Chow motives
      \[   h(X)\cong h(S)(1-r)\oplus \bigoplus_{i=0}^{r-2} h(U)(-i)\ \ \ \hbox{in}\ \MM_{\rm rat}\ .\]
      \end{theorem}

    \begin{proof} This is \cite[Corollary 3.2]{Ji}, which more precisely gives an isomorphism of {\em integral\/} Chow motives. Let us give some details about the isomorphism as constructed in loc. cit.. Let 
    \[ \Gamma:= X\times_U S\ \ \subset\ X\times S \]
  (this is equal to $\PP(\NNN_i)=\HH_s\times_X Z$ in the notation of loc. cit.).  Let
    \[ \Pi_i\ \ \in\ A^\ast(X\times U)\ \ \ \ (i=0, \ldots, r-2) \]
    be correspondences inducing the maps $ (\pi_i)_\ast$ of loc. cit., i.e.  
       \[ (\Pi_i)_\ast= (\pi_i)_\ast :=  (q_{i+1})_\ast \iota_\ast\colon\ \ A^j(X)\ \to\ A^{j-i}(U)\ ,\]
       where $\iota\colon X\hookrightarrow\PP(E)$ is the inclusion morphism, and the $(q_{i+1})_\ast\colon A_\ast(\PP(E))\to A_\ast(U)$ are defined in loc. cit. in terms of the projective bundle formula for $q\colon E\to U$. As indicated in \cite[Corollary 3.2]{Ji} (cf. also \cite[text preceeding Corollary 3.2]{Ji}), there is an isomorphism
       \[   \Bigl( \Gamma, \Pi_0,\Pi_1,\ldots, \Pi_{r-2}\Bigr)\colon\ \ h(X)\ \xrightarrow{\cong}\ h(S)(1-r)\oplus \bigoplus_{i=0}^{r-2} h(U)(-i)\ \ \ \hbox{in}\ \MM_{\rm rat}\ .\]    
    \end{proof}

\begin{remark} In the set-up of Theorem \ref{ji}, a cohomological relation between $X$ and $S$ was established in \cite[Prop. 4.3]{Ko} (cf. also \cite[section 3.7]{IM0}, as well as \cite[Proposition 46]{BFM} for a generalization). A relation on the level of derived categories was established in \cite[Theorem 2.10]{Or} (cf. also \cite[Theorem 2.4]{KKLL} and \cite[Proposition 47]{BFM}).
\end{remark}

   \subsection{The Franchetta property}

 \begin{definition} Let $\XX\to B$ be a smooth projective morphism, where $\XX, B$ are smooth quasi-projective varieties. We say that $\XX\to B$ has the {\em Franchetta property in codimension $j$\/} if the following holds: for every $\Gamma\in A^j(\XX)$ such that the restriction $\Gamma\vert_{X_b}$ is homologically trivial for the very general $b\in B$, the restriction $\Gamma\vert_b$ is zero in $A^j(X_b)$ for all $b\in B$.
 
 We say that $\XX\to B$ has the {\em Franchetta property\/} if $\XX\to B$ has the Franchetta property in codimension $j$ for all $j$.
 \end{definition}
 
 This property is studied in \cite{PSY}, \cite{BL}, \cite{FLV}, \cite{FLV3}.
 
 \begin{definition} Given a family $\XX\to B$ as above, with $X:=X_b$ a fiber, we write
   \[ GDA^j_B(X):=\ima\Bigl( A^j(\XX)\to A^j(X)\Bigr) \]
   for the subgroup of {\em generically defined cycles}. 
  In a context where it is clear to which family we are referring, the index $B$ will often be suppressed from the notation.
  \end{definition}
  
  With this notation, the Franchetta property amounts to saying that $GDA^\ast_B(X)$ injects into cohomology, under the cycle class map.

    \subsection{The CK property}
  
  \begin{definition} Let $M$ be a smooth quasi-projective variety. We say that $M$ has {\em the CK property\/} if for any smooth quasi-projective variety $Z$, the natural map
    \[ A^\ast(M)\otimes A^\ast(Z)\ \to\ A^\ast(M\times Z) \]
    is surjective.
    
    (NB: ``CK'' stands for ``Chow--K\"unneth''.)
     \end{definition}
     
     \begin{lemma}\label{lemCK} Let $\bar{M}$ be a smooth projective variety.
     
     \noindent
     (\rom1) $\bar{M}$ has the CK property if and only if $\bar{M}$ has trivial Chow groups.
     
     \noindent
     (\rom2) if $\bar{M}$ has the CK property, any Zariski open $M\subset \bar{M}$ has the CK property.
     \end{lemma}
  
  \begin{proof}
  
  \noindent
  (\rom1) This is well-known. Assume $\bar{M}$ has the CK property. Then the diagonal $\Delta_{\bar{M}}$ is {\em completely decomposable\/}, i.e. one has
    \begin{equation}\label{ddec} \Delta_{\bar{M}}= V_1\times W_1 + \cdots + V_r\times W_r\ \ \ \hbox{in}\ A^{\dim \bar{M}}(\bar{M}\times\bar{M})\ ,\end{equation}
    where $V_j, W_j\subset\bar{M}$ are subvarieties and $\dim V_j+\dim W_j=\dim\bar{M}$. Looking at the action of this decomposition on Chow groups, one finds that
    $A^\ast_{hom}(X)=0$. 
    
    Conversely, assume $\bar{M}$ has trivial Chow groups. The Bloch--Srinivas argument \cite{BS} then gives a decomposition \eqref{ddec}, which means that the motive of $\bar{M}$ is a sum of trivial motives:
      \[ h(\bar{M})=\oplus \one(\ast)\ \ \ \hbox{in}\ \MM_{\rm rat}\ .\] 
      Given any smooth projective variety $\bar{Z}$, it follows that
      \[  h(\bar{M}\times\bar{Z})=\oplus h(\bar{Z})(\ast)\ \ \ \hbox{in}\ \MM_{\rm rat}\ , \]
      and so
          \[ A^\ast(\bar{M})\otimes A^\ast(\bar{Z})\ \xrightarrow{\cong}\ A^\ast(\bar{M}\times\bar{Z})\ .\]  
          Given a smooth quasi-projective variety $Z$, let $\bar{Z}\supset Z$ be a smooth compactification. The commutative diagram
          \[ \begin{array} [c]{ccc} 
    A^\ast(\bar{M})\otimes A^\ast(\bar{Z})& \xrightarrow{}& A^\ast(\bar{M}\times\bar{Z})\\
    &&\\
    \downarrow && \downarrow\\
    &&\\
      A^\ast(\bar{M})\otimes A^\ast(Z)& \xrightarrow{}& A^\ast(\bar{M}\times Z)\\   
      \end{array} \]
      (where vertical arrows are surjections) implies that the lower horizontal arrow is surjective, i.e. $\bar{M}$ has the CK property.       
          
  \noindent
  (\rom2) This is immediate in view of the commutative diagram
         \[ \begin{array} [c]{ccc} 
    A^\ast(\bar{M})\otimes A^\ast(Z)& \xrightarrow{}& A^\ast(\bar{M}\times Z)\\
    &&\\
    \downarrow && \downarrow\\
    &&\\
      A^\ast(M)\otimes A^\ast(Z)& \xrightarrow{}& \ A^\ast(M\times Z)\ .\\   
      \end{array} \]
  
  \end{proof}
  
  \begin{example}\label{exPf}
  Here is the main example we have in mind: Let $V$ be a vector space of dimension $n$, and consider the {\em Pfaffian\/} of degenerate skew forms
  \[ \Pf:= \Bigl\{ \omega\in \PP(\wedge^2 V^\vee)\, \Big\vert\,  \hbox{rank} \, \omega < r_{max} \Bigr\}\ \ \ \subset\  \PP(\wedge^2 V^\vee)\ \]
  (where $r_{max}=n$ if $n$ is even, and $r_{max}=n-1$ if $n$ is odd). We claim that the non-singular locus
    \[ \Pf^\circ:= \Pf\setminus \Sing(\Pf)\ \]
  has the CK property.
  
  To see this, we consider the variety
    \[   \wt{\Pf}:= \Bigl\{ (\omega,K)\in  \Pf\times \Gr(s,n)   \, \Big\vert\,  K\subset\ker\omega  \Bigr\} \ \ \ \subset\ \Pf\times \Gr(s,n)   \ , \]
    where $s=2$ if $n$ is even, and $s=3$ if $n$ is odd.
    The projection $\wt{\Pf}\to\Gr(2,n)$ is a projective bundle (and so $\wt{\Pf}$ is smooth), and the projection $\wt{\Pf}\to\Pf$ is an isomorphism over the non-singular locus (and so
    $\wt{\Pf}\to\Pf$ is a resolution of singularities). Because $\wt{\Pf}$ (being a projective bundle over a Grassmannian) has trivial Chow groups, the claim follows from Lemma \ref{lemCK}.    
  \end{example}

  \subsection{A Franchetta-type result}

  \begin{proposition}\label{Frtype} Let $M=\bar{M}\setminus N$, where $M$ is a projective variety and $N\subset M$ is closed. Assume $M$ is smooth and has the CK property. Let $L_1,\ldots,L_r\to \bar{M}$ be very ample line bundles, and let
  $\XX\to B$ be the universal family of smooth complete intersections of type 
    \[ X=\bar{M}\cap H_1\cap\cdots\cap H_r\ ,\ \ \  H_j\in\vert L_j\vert\ .\]
  Assume all fibers $X=X_b$ are disjoint from $N$ and have $H^{\dim X}_{tr}(X,\QQ)\not=0$.
  There is an inclusion
    \[ \ker \Bigl( GDA^{\dim X}_B(X\times X)\to H^{2\dim X}(X\times X,\QQ)\Bigr)\ \ \subset\ \Bigl\langle (p_1)^\ast GDA^\ast_B(X), (p_2)^\ast GDA^\ast_B(X)  \Bigr\rangle\ .\]
   \end{proposition}
   
   \begin{proof} In case $M=\bar{M}$ is projective, this is essentially Voisin's ``spread'' result \cite[Proposition 1.6]{V1} (cf. also \cite[Proposition 5.1]{LNP} for a reformulation of Voisin's result). We give a proof which is somewhat different from \cite{V1}. Let $\bar{B}:=\PP H^0(\bar{M},L_1\oplus\cdots\oplus L_r)$ (so $B\subset \bar{B}$ is a Zariski open), and let us consider the projection
   \[ \pi\colon\ \  \XX\times_{\bar{B}} \XX\ \to\ Y\times Y\ .\]
   Using the very ampleness assumption, one finds that $\pi$ is a $\PP^s$-bundle over $(\bar{M}\times \bar{M})\setminus \Delta_{\bar{M}}$, and a $\PP^t$-bundle over $\Delta_{\bar{M}}$.
   That is, the morphism $\pi$ is what is termed a {\em stratified projective bundle\/} in \cite{FLV}. As such, \cite[Proposition 5.2]{FLV} implies the equality
      \begin{equation}\label{stra} GDA^\ast_B(X\times X)= \ima\Bigl( A^\ast(M\times M)\to A^\ast(X\times X)\Bigr) +  \Delta_\ast GDA^\ast_B(X)\ ,\end{equation}
      where $\Delta\colon X\to X\times X$ is the inclusion along the diagonal. As $M$ has the CK property, $A^\ast(M\times M)$ is generated by $A^\ast(M)\otimes A^\ast(M)$. 
      Base-point freeness of the $L_j$ implies that 
        \[  GDA^\ast_B(X)=\ima\bigl( A^\ast(M)\to A^\ast(X)\bigr)\ .\]
       The equality \eqref{stra} thus reduces to
      \[ GDA^\ast_B(X\times X)=\Bigl\langle (p_1)^\ast GDA^\ast_B(X), (p_2)^\ast GDA^\ast_B(X), \Delta_X\Bigr\rangle\ \]   
      (where $p_1, p_2$ denote the projection from $S\times S$ to first resp. second factor). The assumption that $X$ has non-zero transcendental cohomology
      implies that the class of $\Delta_X$ is not decomposable in cohomology. It follows that
      \[ \begin{split}  \ima \Bigl( GDA^{\dim X}_B(X\times X)\to H^{2\dim X}(X\times X,\QQ)\Bigr) =&\\
       \ima\Bigl(  \Dec^{\dim X}(X\times X)\to H^{2\dim X}(X\times X,\QQ)\Bigr)& \oplus \QQ[\Delta_X]\ ,\\
       \end{split}\]
      where we use the shorthand 
       \[ \Dec^j(X\times X):= \Bigl\langle (p_1)^\ast GDA^\ast_B(X), (p_2)^\ast GDA^\ast_B(X)\Bigr\rangle\cap A^j(X\times X) \ \]     
       for the {\em decomposable cycles\/}. 
       We now see that if $\Gamma\in GDA^{\dim X}(X\times X)$ is homologically trivial, then $\Gamma$ does not involve the diagonal and so $\Gamma\in \Dec^{\dim X}(X\times X)$.
       This proves the proposition.
         \end{proof}
  
  \begin{remark} Proposition \ref{Frtype} has the following consequence: if the family $\XX\to B$ has the Franchetta property, then $\XX\times_B \XX\to B$ has the Franchetta property in codimension $\dim X$.
  
  Note that the condition $H^{\dim X}_{tr}(X,\QQ)\not=0$ is necessary to ensure that the class of the diagonal in cohomology is not decomposable. This condition is ``not very restrictive'', cf. \cite[Remark 0.8]{V1} where a similar condition is discussed.
    \end{remark}

 \subsection{A Hilbert schemes argument}
 
  \begin{proposition}[Voisin \cite{V0}, \cite{V1}]\label{spread} Let $\XX$, $\YY$ and $\Zz$ be families over $B$, and assume the morphisms to $B$ are smooth projective and the total spaces are smooth quasi-projective. Let
    \[  \Gamma\in\ \  A^i(\XX\times_B \Zz)\    \]
 be a relative correspondence, with the property that for any $b\in B$ there exist correspondences $\Lambda_b\in A^{\ast}(Y_b\times Z_b)$ and $\Psi_b\in A^\ast(X_b\times Y_b)$ such that
             \[  \Gamma\vert_b= \Lambda_b\circ  \Psi_b\ \ \hbox{in}\ H^{2i}(X_b\times Z_b,\QQ)\ .\]
      Then there exist relative correspondences
      \[ \Lambda\ \ \in A^{\ast}(\YY\times_B \Zz) \ ,\ \ \  \Psi\ \ \in\ A^\ast(\XX\times_B \YY)\   \]
      with the property that for any $b\in B$
      \[ \Gamma\vert_b=\Lambda\vert_b \circ \Psi\vert_b\ \ \hbox{in}\ H^{2i}(X_b\times Z_b,\QQ)\ .\]
   \end{proposition}
   
   \begin{proof} The statement is different, but this is really the same Hilbert schemes argument as \cite[Proposition 3.7]{V0}, \cite[Proposition 4.25]{Vo} (cf. also \cite[Proposition 2.10]{Lfam} for a similar statement). The point is that the data of 
   all the
   $(b,\Lambda_b, \Psi_b)$ that are solutions to the splitting problem
   \[   \Gamma\vert_b= \Lambda_b\circ \Psi_b\ \ \hbox{in}\ H^{2i}(X_b\times Z_b,\QQ)\ \]
   can be encoded by a countable number of algebraic varieties $p_j\colon M_j\to B$, with universal objects 
     \[   U_j= (\Lambda_j,\Psi_j)\ ,\ \ \ \ \Lambda_j\subset \YY\times_{M_j}\Zz\ , \Psi_j\subset\XX\times_{M_j}\YY\  ,\] 
     with the property that 
   for $m\in M_j$ and $b=p_j(m)\in B$, we have
     \[  (U_j)\vert_{m}=(\Lambda_b,\Psi_b)\ \ \hbox{in}\ H^{\ast}(Y_b\times Z_b)\oplus H^\ast(X_b\times Y_b)\ .\]
     By assumption, the union of the $M_j$ dominates $B$. Since there is a countable number of $M_j$, one of the $M_j$ (say $M_0$) must dominate $B$. Taking hyperplane sections, we may assume $M_0\to B$ is generically finite (say of degree $d$). Projecting $\Lambda_0$ to $\YY\times_B \Zz$ (resp. projecting $\Psi_0$ to $\XX\times_B \YY$) and dividing by $d$, we have obtained $\Lambda$ (resp. $\Psi$) as requested.
   \end{proof}   

This is a variant of Proposition \ref{spread}:

 \begin{proposition}\label{spread2} Let $\XX$, $\YY$ be families over $B$, and assume the morphisms to $B$ are smooth projective and the total spaces are smooth quasi-projective. Let
    \[  \Gamma_\XX\in\ \  A^i(\XX\times_B \XX)\  ,\ \ \ \Gamma_\YY\in\ \ A^j(\YY\times_B \YY)    \]
 be relative correspondences, with the property that for any $b\in B$ there exist correspondences $\Lambda_b\in A^{\ast}(Y_b\times X_b)$ and $\Psi_b\in A^\ast(X_b\times Y_b)$ such that
             \[  \begin{split} \Gamma_\XX\vert_b&= \Lambda_b\circ  \Psi_b\ \ \hbox{in}\ H^{2i}(X_b\times X_b,\QQ)\ ,\\
                                          \Gamma_\YY\vert_b&= \Psi_b\circ  \Lambda_b\ \ \hbox{in}\ H^{2j}(Y_b\times Y_b,\QQ)\ ,\\    
                                          \end{split}         \]
      Then there exist relative correspondences
      \[ \Lambda\ \ \in A^{\ast}(\YY\times_B \XX) \ ,\ \ \  \Psi\ \ \in\ A^\ast(\XX\times_B \YY)\   \]
      with the property that for any $b\in B$
        \[  \begin{split} \Gamma_\XX\vert_b&= \Lambda\vert_b\circ  \Psi\vert_b\ \ \hbox{in}\ H^{2i}(X_b\times X_b,\QQ)\ ,\\
                                          \Gamma_\YY\vert_b&= \Psi\vert_b\circ  \Lambda\vert_b\ \ \hbox{in}\ H^{2j}(Y_b\times Y_b,\QQ)\ ,\\    
                                          \end{split}         \]      
       \end{proposition}
 
 \begin{proof} As in the proof of Proposition \ref{spread}, the data of all $(b,\Lambda_b, \Psi_b)$ are encoded by a countable number of algebraic varieties carrying universal objects. The same argument then applies.
 \end{proof}

 \subsection{Transcendental motive and variable motive}
 
 \begin{proposition}\label{trans} Let $X$ be a smooth projective variety. There exists a splitting
   \[ h(X)=t(X) \oplus h_{alg}(X)\ \ \ \hbox{in}\ \MM_{\rm rat}\ ,\]
   such that
    \[ H^\ast(t(X),\QQ)= H^\ast_{tr}(X,\QQ) \]
    (where $H^\ast_{tr}(X,\QQ)$ is defined as the orthogonal complement of the algebraic part of cohomology $H^\ast_{alg}(X,\QQ)$), and
    \[ A^\ast(t(X))=A^\ast_{hom}(X)\ .\]
      \end{proposition}
      
      \begin{proof} This is a standard construction, cf. for instance \cite{43} where projectors on $H^\ast_{alg}(X,\QQ)$ are constructed.      
      \end{proof}
 
 The inconvenience of the decomposition of Proposition \ref{trans} is that $t(X)$ is not canonically defined (the construction depends on choices). Also, if $X$ varies in a family it is not clear if the decomposition is generically defined. For complete intersections $X\subset M$, there is another decomposition which does {\em not\/} present these inconveniences. This decomposition
 is a motivic incarnation of the so-called ``variable cohomology'', i.e. the cohomology of $X$ not coming from the ambient space $M$ (for more on this notion, and the relation with primitive cohomology, cf. \cite[Section 4.3.4]{Vo} or \cite{Pet}).

 \begin{proposition}\label{var} Let $M$ be a smooth projective variety with trivial Chow groups, and let $X\subset M$ be a smooth complete intersection (defined by very ample line bundles on $M$).
 There is a canonical splitting
   \[ h(X) = h_{\rm var}(X)\oplus h_{\rm fix}(X)\ \ \ \hbox{in}\ \MM_{\rm rat}\ ,\]
   such that
      \[  H^\ast(h_{\rm var}(X),\QQ) = H^{\dim X}_{\rm var}(X,\QQ):= \coker\Bigl( H^{\dim X}(M,\QQ)\to  H^{\dim X}(X,\QQ)\Bigr)\ .\]
    
    Moreover, there exists $t(X)$ as in Proposition \ref{trans} such that there is an inclusion as submotive
    \[  t(X)  \ \hookrightarrow\ h_{\rm var}(X) \ \ \ \hbox{in}\ \MM_{\rm rat}\ .\]      
        \end{proposition} 
        
        \begin{proof} This is a standard construction, cf. for instance \cite{Pet} (where the hypotheses on $M$ are less stringent).
        
        For the ``moreover'' part, one can start from $h_{\rm var}(X)$ and take out trivial motives generating the algebraic part of $H^\ast(h_{\rm var}(X),\QQ)$; the difference $t(X):= h_{\rm var}(X)-\oplus \one(\ast)$ is as requested.      
                  \end{proof}

 \section{The Pfaffian--Grassmannian equivalence}
 
 \subsection{Set-up}
 
 \begin{notation}\label{not} Given $V$ a vector space of dimension $n$, let
  \[   \Gr(2,V)\ \ \ \subset\ \PP(\wedge^2 V)   \]
  denote the Grassmannian of 2-dimensional linear subspaces in $V$ in its Pl\"ucker embedding.
  The projective dual to $ \Gr(2,V)$ is the {\em Pfaffian\/} of degenerate skew forms
  \[ \Pf:= \Bigl\{ \omega\in \PP(\wedge^2 V^\vee)\, \Big\vert\,  \hbox{rank} \, \omega < n \Bigr\}\ \ \ \subset\  \PP(\wedge^2 V^\vee)\ .\]
 Given a linear subspace $U\subset \wedge^2 V$ of codimension $k$, we define varieties by intersecting
    \[  \begin{split}   X&=X_U:=  \Gr(2,V)\cap   \PP(U) \ \ \ \subset\ \PP(\wedge^2 V) \ , \\
                              Y&=Y_U:= \Pf\,\cap\, \PP(U^\perp)\ \ \ \subset\ \PP(\wedge^2 V^\vee) \ .
                              \end{split}  \]  
  For $U$ generic and $k$ small enough (the precise condition is that $k\le 6$ when $n$ is even, and $k\le 10$ when $n$ is odd, so that $Y$ avoids the singular locus of $\Pf$), the intersections $X$ and $Y$ are smooth and dimensionally transverse, of dimension
   \[  \begin{split}  \dim X&=2(n-2)-k\ ,\\
                             \dim Y&=\begin{cases} k-2 &\hbox{if\ $n$\ even},\\
                                                                  k-4 &\hbox{if\ $n$\ odd}.\\    
                                                                  \end{cases} \\
                                                                  \end{split} \]
  \end{notation}

 \subsection{Relation in the Grothendieck ring}
 \label{ss:l}

\begin{lemma}\label{lemK0} Assume $n\ge 4$.
There is a relation in $K_0(Var_\C)$
     \[ [\Gr(2,n)]= \begin{cases}  [\PP^{n-2}]\cdot {\displaystyle\sum_{k=0}^{(n-2)/2}} \LLL^{2k}\ &\ \hbox{if\ $n$\ even}\ ,\\
                                             [\PP^{n-1}]\cdot {\displaystyle\sum_{k=0}^{(n-3)/2}} \LLL^{2k}\ &\ \hbox{if\ $n$\ odd}\ .\\
                                             \end{cases}\]
   Moreover, let $H(2,n)\subset\Gr(2,n)$ denote a smooth Pl\"ucker hyperplane section. There is a relation in $K_0(Var_\C)$
     \[ [H(2,n)]= \begin{cases}  [\PP^{n-3}]\cdot {\displaystyle\sum_{k=0}^{(n-2)/2}} \LLL^{2k}\ &\ \hbox{if\ $n$\ even}\ ,\\
                                             [\PP^{n-2}]\cdot {\displaystyle\sum_{k=0}^{(n-3)/2}} \LLL^{2k}\ &\ \hbox{if\ $n$\ odd}\ .\\
                                             \end{cases}\]
    \end{lemma}                                         

\begin{proof} The formula for $[\Gr(2,n)]$ is \cite[Section 2]{Mar}. As for the second formula, it is well-known that there is no variable cohomology:
 \[ H^{2(n-2)-1}(H(2,n),\QQ)=0\ \] 
 \cite[Proposition 2.3]{Don}.
Since the cohomology in degree $<2(n-2)-1$ is isomorphic to that of $\Gr(2,n)$ (weak Lefschetz), this gives the formula for $[H(2,n)]$.
\end{proof}

\begin{theorem}\label{K0} Given $n,k\in\NN$ such that $n$ is even and $k\le 6$, or $n$ is odd and $k\le 10$, let 
  \[  X\ \subset \ \Gr(2,V)\ ,\ \ \ Y\ \subset\  \Pf \]
  be smooth dimensionally transverse intersections as in Notation \ref{not}.
There is a relation in the Grothendieck ring of varieties
  \[  [X]\cdot\LLL^{k-1} + [\PP^{k-2}]\cdot [\Gr(2,n)]  = [Y]\cdot \LLL^s + [\PP^{k-1}]\cdot [H(2,n)]
  \ \ \ \hbox{in}\ K_0(Var_\C)\ ,\]
  where 
    \[ s=\begin{cases} n-1 &\ \hbox{if\ $n$\ odd}\ ,\\
                            n-2 &\ \hbox{if\ $n$\ even}\ .\\
                            \end{cases}\]
                            \end{theorem}
                            
     \begin{proof} (This proof is inspired by \cite{Mar}, where the case $(n,k)=(7,7)$ is done.)
     
     Let us consider 
      \[ Q:= \bigl\{ (T,\C\omega)\ \in\ \Gr(2,V)\times \PP(U^\perp)\  \big\vert\   \omega\vert_T=0\bigr\}\ \ \subset\  \Gr(2,V)\times \PP(U^\perp)\  , \]
                  the so-called {\em Cayley hypersurface\/}.  
     There is a diagram
   \begin{equation}\label{diag} \begin{array}[c]{ccccccccc}   && Q_X & \hookrightarrow & Q & \hookleftarrow & Q_Y && \\
                &&&&&&&&\\
                   &{}^{\scriptstyle } \swarrow \ \ && {}^{\scriptstyle p} \swarrow \ \ \ & & \ \ \ \searrow {}^{\scriptstyle q} & & \ \ \searrow {}^{\scriptstyle } & \\
                   &&&&&&&&\\
                   X & \hookrightarrow & \Gr(2,V) &  &  & & \PP(U^\perp) & \hookleftarrow & Y\\
                   \end{array}\end{equation}
                     Here, the morphisms $p$ and $q$ are induced by the natural projections, and the closed subvarieties $Q_X, Q_Y\subset Q$ are defined as $p^{-1}(X)$ resp. $q^{-1}(Y)$. 

The restriction of $p$ to $Q\setminus Q_X$ is trivial with fibre $Q_u\cong\PP^{k-2}$, while the restriction of $p$ to $Q_X$ is Zariski locally trivial with fibre $Q_{X,x}\cong\PP^{k-1}$. This allows us to relate $Q$ and $X$ in the Grothendieck ring:
 \begin{equation}\label{easyside}  \begin{split}       [Q]&=  [Q\setminus Q_X] + [Q_X]\\
                                      &= [\Gr(2,n)\setminus X]\cdot [\PP^{k-2}] + [X]\cdot [\PP^{k-1}]\\
                                      &=  \Bigl( [\Gr(2,n)]-[ X]\Bigr)\cdot \sum_{j=0}^{k-2} \LLL^j + [X]\cdot     \sum_{j=0}^{k-1}   \LLL^j\\
                                      &= [\Gr(2,n)] \cdot [\PP^{k-2}] + [X] \cdot \LLL^{k-1}\ \ \ \ \ \ \hbox{in}\ K_0(Var_\C)\ .
                                      \end{split}\end{equation}
This gives us the left-hand side of the relation of Theorem \ref{K0}.

For the right-hand side of the relation, we express $[Q]$ in terms of $[Y]$ by exploiting the right-hand side of diagram \eqref{diag}. This is more involved (for one thing, the morphisms on the right-hand side are not known to be Zariski locally trivial fibrations), but can be done with some patience. The first thing to note is that the morphism $q$ has only 2 types of fibers:

\begin{lemma}\label{pt} The morphism $Q\setminus Q_Y\to  \PP(U^\perp)\setminus Y$ (obtained by restricting $q$) is piecewise trivial (in the sense of \cite[Section 4.2]{Seb}) with constant fiber $F_1$.

Likewise, the morphism $ Q_Y\to  Y$ (obtained by restricting $q$) is piecewise trivial with constant fiber $F_2$.
\end{lemma}        

 \begin{proof} The argument of \cite[Lemma 3.3]{Mar} for the case $(n,k)=(7,7)$ extends. Let us give a sketch for the morphism $ Q_Y\to  Y$.   
 In view of \cite[Theorem 4.2.3]{Seb}, it suffices to check that there exists $F_2$ such that for all $y\in Y$ the fiber $q^{-1}(y)$ is a $\C(x)$-scheme isomorphic to $F_2 \times_\C \C(x)$. 
 To check this, one observes that a skew form of rank $2r$ (where $2r=n-2$ if $n$ is even, and $2r=n-3$ if $n$ is odd) with coefficients in a field $K\supset\C$ is congruent to the skew form
   \[         \begin{bmatrix}  0_r & I_r &    \\
                                              -I_r & 0_r &  \\
                                               & & 0_{n-2r} \\
                                               \end{bmatrix}\]
                                               (where $0_r$ and $I_r$ denote the $r\times r$ zero-matrix resp. identity matrix, and all empty entries are zero), with a base change having coefficients in $K$.
  \end{proof}      
  
  The next step consists in finding expressions for the classes of the fibers $F_1$ and $F_2$ in the Grothendieck ring:
  
  \begin{lemma}\label{F} One has
     \[ \begin{split}  [F_1]&=[H(2,n)]\ ,\\
                       [F_2]&= [H(2,n)] + \LLL^s
     \ \ \ \ \ \hbox{in}\ K_0(Var_\C)\ ,
     \end{split} \]
     where
       \[  s=\begin{cases} n-1 &\ \hbox{if\ $n$\ odd}\ ,\\
                            n-2 &\ \hbox{if\ $n$\ even}\ .\\
                            \end{cases}\]
     \end{lemma}     
     
     \begin{proof} For $F_1$, we note that the open $\PP(U^\perp)\setminus Y$ is exactly the locus where $q$ has smooth fibers, and (since the Cayley hypersurface $Q\subset  \Gr(2,n)\times \PP(U^\perp)$ is a $(1,1)$-divisor) a fiber of $q$ is a Pl\"ucker hyperplane section of $\Gr(2,n)$. (Alternatively, one can do a cut-and-paste argument similar to \cite[Lemma 3.5]{Mar}, which is the case $(n,k)=(7,7)$.)     
     
     For $F_2$, we can apply the stratification argument of \cite[Lemma 3.4]{Mar} (which is the case $(n,k)=(7,7)$); a little computation gives the formula. Alternatively, we can reason as in \cite[Proposition 2.6]{Man}: a singular Pl\"ucker hyperplane section $F_2$ of $\Gr(2,n)$ is isomorphic to a Schubert divisor, and so $F_2$ has a cell decomposition given by all the cells of 
$\Gr(2,n)$ minus the big cell. It follows that
  \[ [F_2]= [\Gr(2,n)] - \LLL^{2n-4}\ \ \ \hbox{in}\ K_0(Var_\C)\ ,\]
  and the formula follows from Lemma \ref{lemK0}.
      \end{proof}       
     
     Armed with Lemmata \ref{pt} and \ref{F}, we can relate $Q$ and $Y$ in the Grothendieck group:
       \begin{equation}\label{hardside} \begin{split}   [Q] &=  [Q\setminus Q_Y] + [Q_Y]\\
                                                                 &= [\PP^{k-1}\setminus Y]\cdot [F_1]  +  [Y]\cdot [F_2] \\
                                                                 &= \Bigl(   [\PP^{k-1}] - [Y]\Bigr)\cdot [H(2,n)]  + [Y]\cdot \Bigl(  [H(2,n)] + \LLL^s\Bigr)\\
                                                                 &=      [\PP^{k-1}] \cdot   [H(2,n)]  +    [Y]\cdot  \LLL^s   \ \ \ \ \ \ \hbox{in}\ K_0(Var_\C)\ .
                                      \end{split}\end{equation}
  (Here we have used the general fact that if $M\to N$ is a piecewise trivial morphism with constant fiber $F$, one has equality
        \[ [M]= [N]\cdot [F] \ \ \ \hbox{in}\ K_0(Var_\C)\ ,\]
        which is \cite[Lemma 4.2.2]{Seb}.)

Theorem \ref{K0} now follows from writing
   \[    [Q]=[Q]\ \ \    \hbox{in}\ K_0(Var_\C)\ ,\]
   and applying equalities \eqref{easyside} and \eqref{hardside}.
                                \end{proof}

\begin{corollary}\label{corK0} Let $X, Y$ be as in Theorem \ref{K0}, and assume $n=k$ is odd. Then $X$ and $Y$ are Calabi--Yau varieties of dimension $n-4$ and 
   \[  ([X]-[Y])\cdot \LLL^{k-1}=0\ \ \ \hbox{in}\ K_0(Var_\C)\ \]
   (i.e., $X$ and $Y$ are L-equivalent).
   \end{corollary}
   
   \begin{proof} This follows from Theorem \ref{K0}, combined with the fact that
     \[ [ \PP^{n-1}]\cdot [H(2,n)] =   [ \PP^{n-2}]\cdot    [\Gr(2,n)]  \ \ \ \hbox{in}\ K_0(Var_\C)\ . \]
     (this equality follows from Lemma \ref{lemK0}).     
   \end{proof}

\begin{remark} In case $n=k=7$, $X$ and $Y$ are the famous Calabi--Yau threefolds of Borisov \cite{Bor}. In this case, the relation of Corollary \ref{corK0} was proven by 
Martin \cite{Mar}, improving on work of Borisov. (For the fact that in general $[X]\not=[Y]$ in the Grothendieck ring, so that $\LLL$ is a zero-divisor, this follows from the fact that $X$ and $Y$ are not stably birational, cf. \cite[Proof of Theorem 2.12]{Bor}).

For $n=k=9$, $X$ and $Y$ are Calabi--Yau fivefolds that are derived equivalent (cf. Theorem \ref{hpd}(\rom3)). For this case, the L-equivalence of Corollary \ref{K0} appears to be new.
\end{remark}

\subsection{Relation in cohomology}
\label{ss:coh}

The relation in the Grothendieck ring (Theorem \ref{K0}) has consequences for the cohomology of the Pfaffian--Grassmannian equivalence.
We first deduce a numerical statement:

\begin{lemma}\label{corK02} Assume $n$ is even and $k\in\{2,4\}$, or $n$ is odd and $k\in\{2,4, 6\}$. Let
  \[  X\ \subset \ \Gr(2,V)\ ,\ \ \ Y\ \subset\  \Pf \]
  be smooth dimensionally transverse intersections as in Notation \ref{not}.
Then
  \[ \dim H^{\dim X}_{var}(X,\QQ) = \dim H^{\dim Y}_{}(Y,\QQ)-1\ .\]
\end{lemma}

\begin{proof} Let us write $k=2 k^\prime$. Since there is a functor from $K_0(Var_\C)$ to the category of Hodge structures, the relation in $K_0(Var_\C)$ of Theorem \ref{K0} gives an equality of Betti numbers
  \[  b_{\dim X}(X) +\sum_{j= n-1-k^\prime} ^{n-3+k^\prime} b_{2j}(\Gr(2,n)) = b_{\dim Y}(Y) + \sum_{j= n-2-k^\prime}^{n-3+k^\prime} b_{2j}(H(2,n)) \ \]
  (where for any variety $M$ we write $b_j(M):=\dim H^j(M,\QQ)$).
  
Using the equalities $b_j(\Gr(2,n))=b_j(H(2,n))$ for $j<2n-4$ and $b_j(H(2,n))=b_{j+2}(\Gr(2,n))$ for $j\ge 2n-4$ (weak Lefschetz), this simplifies to
\[     b_{\dim X}(X) + b_{2n-4}(\Gr(2,n)) =  b_{\dim Y}(Y) + b_{2n+2k^\prime-6}(H(2,n))+  b_{2n-2k^\prime-4}(H(2,n))\ .\]
Expressing everything in terms of Betti numbers of $\Gr(2,n)$ (weak Lefschetz), this becomes
  \[   b_{\dim X}(X) + b_{2n-4}(\Gr(2,n)) =  b_{\dim Y}(Y) +  2 b_{2n-4-2k^\prime}(\Gr(2,n))\ .\]
Using the equality
  \[  b_{2n-4}(\Gr(2,n))-  b_{2n-4-2k^\prime}(\Gr(2,n))  = \begin{cases}  \lceil {k^\prime\over 2}\rceil   & \hbox{if\ $n$\ is\ even},\\
                                                                                                                    \lfloor {k^\prime\over 2}\rfloor   & \hbox{if\ $n$\ is\ odd}\\
                                                                                                                    \end{cases}\]
                                                                                                                    (which can be deduced from Lemma \ref{lemK0}), one finds that
    \[      b_{\dim X}(X) + 1 =  b_{\dim Y}(Y) +   b_{2n-4-2k^\prime}(\Gr(2,n))\ ,\]
    provided either $n$ is even and $k\in\{2,4\}$, or $n$ is odd and $k\in\{2,4, 6\}$. Since 
    \[      \dim H^{\dim X}_{var}(X,\QQ)=     b_{\dim X}(X) -    b_{2n-4-2k^\prime}(\Gr(2,n)) \]
    this proves the lemma.                                                                                                    
\end{proof}

Restricting to the {\em transcendental cohomology\/}, one obtains a stronger statement: there is a {\em correspondence-induced\/} isomorphism of transcendental cohomology:

\begin{proposition}\label{propcoh} Assume either $k\le 6$, or  $(n,k)=(7,7)$. Let
  \[  X\ \subset \ \Gr(2,V)\ ,\ \ \ Y\ \subset\  \Pf \]
  be smooth dimensionally transverse intersections as in Notation \ref{not}.
There is an isomorphism in cohomology
    \[ H^{\dim X}_{tr}(X,\QQ)\cong H^{\dim Y}_{tr}(Y,\QQ)\ .\]
    (Here $H^\ast_{tr}(X,\QQ)=H^{\dim X}_{tr}(X,\QQ)$ is defined as the orthogonal complement of the algebraic part of cohomology. In particular, $H^{\dim X}_{tr}(X,\QQ)=H^{\dim X}_{}(X,\QQ)$ if $X$ is of odd dimension.)
    
    More precisely, there is an isomorphism of homological motives
    \[ t(X)\ \xrightarrow{\cong}\ t(Y)(-m)\ \ \ \hbox{in}\ \MM_{\rm hom}\ ,\]
   where $m:={1\over 2}(\dim X-\dim Y)$.
  \end{proposition}
 
 \begin{proof} Theorem \ref{K0} gives a relation in the Grothendieck ring
    \[  [X]\cdot\LLL^{k-1} +  \sum_j \LLL^{n_j}   = [Y]\cdot \LLL^s + \sum_j \LLL^{m_j}
  \ \ \ \hbox{in}\ K_0(Var_\C)\ ,\] 
   where 
    \[ s=\begin{cases} n-1 &\ \hbox{if\ $n$\ odd}\ ,\\
                            n-2 &\ \hbox{if\ $n$\ even}\ .\\
                            \end{cases}\]
  Since there is a functor $K_0(Var_\C)\to K_0(\MM_{\rm num})$ sending $[M]$ to $h(M)$ (for any smooth projective variety $M$) and $\LLL^j$ to $\one(-j)$, this gives a relation
   \[ h(X) +\bigoplus\one(\ast) = h(Y)(-m)+  \bigoplus\one(\ast) \ \ \ \hbox{in}\    K_0(\MM_{\rm num})\ \]   
   (where we have used the equality $m=s-k+1$).    
   Because the category $\MM_{\rm num}$ is semi-simple \cite{Jans}, the above induces an isomorphism of motives
    \[ h(X) \oplus\bigoplus\one(\ast) \cong h(Y)(-m)\oplus  \bigoplus\one(\ast) \ \ \ \hbox{in}\   \MM_{\rm num}\ .\ \]   
   Using that (by construction) $h(X)=t(X)\oplus \oplus\one(\ast)$, this gives an isomorphism
     \[ t(X) \oplus\bigoplus_j\one(r_j) \cong t(Y)(-m)\oplus  \bigoplus_k\one(s_k) \ \ \ \hbox{in}\  \MM_{\rm num}\ . \]    
     
    At this point, we note that $X$ and $Y$ verify the standard conjectures (for $X$ this is just because all complete intersections in Grassmannians verify the standard conjectures; for $Y$ this is clear when $n$ is even [because then $Y\subset\PP^{k-1}$ is a hypersurface], and when $n$ is odd it is true because $\dim Y\le 2 $ or $Y$ is a 3-dimensional Calabi--Yau variety, and the standard conjectures are known for threefolds not of general type \cite{Tan}). In particular, homological and numerical equivalence coincide for all self-powers and products of $X$ and 
 $Y$. It follows that the above relation actually takes place in the subcategory $\MM_{\rm hom}^\circ$ generated by homological motives of varieties satisfying the standard conjectures:
   \[  t(X) \oplus\bigoplus_j\one(r_j) \cong t(Y)(-m)\oplus  \bigoplus_k\one(s_k)\ \ \ \hbox{in}\  \MM^\circ_{\rm hom}\ . \]    
   Taking algebraic cohomology $H^\ast_{alg}()$ on both sides (which excludes the transcendental motives), we find that
   \[ \bigoplus_j\one(r_j) \cong  \bigoplus_k\one(s_k) \ \ \ \hbox{in}\  \MM^\circ_{\rm hom}\ . \]       
   Since $ \MM_{\rm hom}^\circ\subset\MM_{\rm num}$ is semi-simple, it follows that there is also an isomorphism
     \[  t(X) \cong t(Y)(-m)\ \ \ \hbox{in}\  \MM^\circ_{\rm hom}\ . \]     
     Since $\MM^\circ_{\rm hom}$ is a subcategory of $\MM_{\rm hom}$, this proves the result.
 \end{proof}
 
 \begin{remark} Let $X,Y$ be as in Proposition \ref{propcoh}, and assume $n$ is even and $k=3$ (i.e. $X$ has dimension $2n-7$ and $Y$ is a curve). In this case, Donagi has proven in his thesis that the intermediate Jacobian of $X$ is naturally isomorphic to the Jacobian of $Y$ \cite[Theorem 2.5]{Don}.
  \end{remark}

 \subsection{Relation of derived categories}
 \label{ss:der} 
 
 The following is included merely for illustrative purposes; this will {\em not\/} be used in this paper: 
 
\begin{theorem}[Kuznetsov \cite{Kuz0}, Segal--Thomas \cite{ST}, Rennemo--Segal \cite{RS}]\label{hpd} Let   
   \[  X\ \subset \ \Gr(2,V)\ ,\ \ \ Y\ \subset\  \Pf \]
  be smooth dimensionally transverse intersections as in Notation \ref{not}.
  
\noindent
(\rom1) Assume $n=6$ or $n=7$. There exist semi-orthogonal decompositions
    \[   \begin{split} D^b(X)&=\bigl\langle E_1,\ldots, E_r, \Kuz(X)\bigr\rangle\ ,\\
                              D^b(Y)&=\bigl\langle E_1,\ldots, E_s, \Kuz(Y)\bigr\rangle\ \\
                              \end{split}\]
                              where the $E_j$ are exceptional objects,
                              and there is an equivalence of categories 
                              \[ \Kuz(X)\ \cong\ \Kuz(Y)\ .\]
                              
 \noindent
 (\rom2) Assume $n$ is odd and $k\le \min(n,10)$, or $n$ is even and $k\le \min(n/2,6)$. Then                              
  there is an embedding
    \[  D^b(Y)\ \hookrightarrow\ D^b(X) \]
    admitting a right adjoint.
    
 \noindent
 (\rom3) Assume $(n,k)=(9,9)$. Then there is an equivalence of derived categories
   \[ D^b(Y)\ \cong\ D^b(X)\ .\]
    \end{theorem}

\begin{proof} Point (\rom1) is one of the first instances of the famous theory of {\em homological projective duality\/} \cite{Kuz0}.

Point (\rom2) is \cite[Theorem 2.9]{ST}; the argument is motivated by but distinct from HPD, which explains why the statement is weaker than that of (\rom1): it is still an open question whether the orthogonal to $D^b(Y)$ inside $D^b(X)$ is generated by exceptional objects (cf. \cite[Remark 3.8]{ST}).

Point (\rom3) is a special case of \cite[Theorem 1.1]{RS}.
\end{proof}

\begin{remark} Conjecturally, the HPD program applies in full generality, and so Theorem \ref{hpd}(\rom1) should be true for all $(n,k)$. This is \cite[Conjecture 5]{Kuz0} (cf. also \cite[Section 5.3]{Thom} and the introduction of \cite{ST}).
\end{remark}

\subsection{The Noether--Lefschetz condition}

\begin{definition}\label{families} Given $V$ a vector space of dimension $n$, let
   \[ \XX\to B\ ,\ \ \ \YY\to B \]
   denote the universal families of smooth dimensionally transverse intersections of type
    \[  \begin{split}   X&=X_U:=  \Gr(2,V)\cap   \PP(U) \ \ \ \subset\ \PP(\wedge^2 V) \ , \\
                              Y&=Y_U:= \Pf\,\cap\, \PP(U^\perp)\ \ \ \subset\ \PP(\wedge^2 V^\vee) \ ,
                              \end{split}  \]   
                              where $U\subset \wedge^2 V$ is a codimension $k$ linear subspace. (In particular, $k\le 6$ if $n$ is even, and $k\le 10$ if $n$ is odd.)
                              
        We write $B^\circ\subset B$ for the Zariski open over which both $X_b$ and $Y_b$ are smooth dimensionally transverse.                      
     \end{definition}                         
                              
  \begin{definition}\label{NL} We say that the family $\XX\to B$ satisfies the condition (NL) if the following holds: for the very general fiber $X_b$, the inclusion
    \[ H^{\dim X_b}_{var}(X_b,\QQ):=\coker\Bigl( H^{\dim X_b}(\Gr(2,V),\QQ)\to H^{\dim X_b}(X_b,\QQ)    \Bigr)       \ \supset\ H^{\dim X_b}_{tr}(X_b,\QQ)     \]
    is an equality, i.e. all Hodge classes in $H^{\dim X_b}(X_b,\QQ)$ come from $\Gr(2,V)$.
    
    (NB: ``NL'' stands for Noether--Lefschetz.)
    
    Likewise, we say that $\YY\to B$ satisfies (NL) if the following holds: for the very general fiber $Y_b$, the inclusion
    \[ H^{\dim Y_b}_{var}(Y_b,\QQ):=\coker\Bigl( H^{\dim Y_b}(  (\PP(\wedge^2 V^\vee) ,\QQ)\to H^{\dim Y_b}(Y_b,\QQ)    \Bigr)       \ \supset\ H^{\dim Y_b}_{tr}(Y_b,\QQ)     \]
    is an equality.
    
    Given 2 integers $(n,k)$, we say the pair $(n,k)$ satisfies condition (NL) if the family $\XX\to B$ and the family $\YY\to B$ both satisfy (NL).    
     \end{definition}
    
    The condition (NL) is trivially fulfilled in case $k$ is odd. Also, at least for small $k$ it suffices to test condition (NL) on one side of the Pfaffian--Grassmanian equivalence:
    
    \begin{lemma}\label{NLequiv} Assume $n$ is even and $k\le 4$, or $n$ is odd and $k\le 6$.
    The family $\XX\to B$ satisfies condition (NL) if and only if $\YY\to B$ satisfies (NL).
    \end{lemma}
    
    \begin{proof} We may restrict to the common base $B^\circ$. For $b\in B^\circ$ very general, we look at the diagram
      \[ \begin{array}[c]{ccc}     H^{\dim X_b}_{tr}(X_b,\QQ) & \xrightarrow{\cong} &   H^{\dim Y_b}_{tr}(Y_b,\QQ) \\
                                           &&\\
                                           \hookdownarrow&&\hookdownarrow\\
                                           &&\\
                                            H^{\dim X_b}_{\rm var}(X_b,\QQ) &   &  \  \ H^{\dim Y_b}_{\rm var}(Y_b,\QQ)\ , \\
                                            \end{array} \]
                                            where the horizontal arrow is the isomorphism of Proposition \ref{propcoh}, and the vertical arrows are the inclusions. Under the hypothesis on $k$, we know that the two spaces at the bottom of this diagram have the same dimension (Lemma \ref{corK02}). It follows that if one of the vertical arrows is an isomorphism, the other vertical arrow is an isomorphism as well.
       \end{proof}
       
%
      
   Here are some examples:
   
   \begin{lemma}\label{lnl} The following pairs satisfy condition (NL):
   
   \noindent
   (0) $(n,k)$ with $k$ odd;
   
   \noindent
   (\rom1) $(2m,4)$ where $m\ge 4$ and $(2m+1,6)$ where $m\ge 3$;
   
   \noindent
   (\rom2) $(6,6), (7,8)$.
    \end{lemma}   
      
   \begin{proof} 
   (\rom1) The assumptions imply that $Y_b$ is a surface with $H^{2,0}(Y_b)\not=0$ (indeed, the canonical bundle of $Y_b$ is either trivial or ample). For the very general $b$, there is then an isomorphism
     \[  A^1(\Pf^\circ) \cong   A^1(Y_b)\ , \]
     this follows from \cite{Mois}.
     
   \noindent
   (\rom2) The assumptions imply that $X_b$ is a surface with $H^{2,0}(X_b)\not=0$. The argument of (\rom1) then applies to the very general $X_b$.  
    \end{proof}

 \subsection{Relation of Chow motives}
 
We now proceed to prove the main result of this paper. In view of Lemma \ref{lnl}, the theorem stated in the introduction is a special case of the following:
 
 \begin{theorem}\label{main} Given $V$ an $n$-dimensional vector space, let 
     \[  \begin{split}   X&=X_U:=  \Gr(2,V)\cap   \PP(U) \ \ \ \subset\ \PP(\wedge^2 V) \ , \\
                              Y&=Y_U:= \Pf\,\cap\, \PP(U^\perp)\ \ \ \subset\ \PP(\wedge^2 V^\vee) \ .
                              \end{split}  \]  
 be smooth dimensionally transverse intersections, where $U\subset \wedge^2 V$ is a codimension $k$ linear subspace. Assume $k\le 6$ or $(n,k)=(7,7)$.
 Assume also that $(n,k)$ satisfies condition (NL), and that
 $H^\ast_{tr}(X,\QQ)\not=0$. Then there exist $t(X), t(Y)$ as in Proposition \ref{trans} such that there
  is an isomorphism of Chow motives
   \[ t(X)\ \xrightarrow{\cong}\ t(Y)(-m)\ \ \ \hbox{in}\ \MM_{\rm rat}\ ,\]
   where $m={1\over 2}(\dim X -\dim Y)$.    
   
   In particular, there are isomorphisms
     \[ A^j_{hom}(X)\cong  A^{j-m}_{hom}(Y)\ \ \ \forall j\ .\]                         
  \end{theorem}
  
  \begin{proof} Thanks to Proposition \ref{propcoh}, we know that there is an isomorphism of homological motives
    \[ \Psi\colon\ \  t(X)\ \xrightarrow{\cong}\ t(Y)(-m)\ \ \ \hbox{in}\ \MM_{\rm hom}\ .\]
  Let us now consider things family-wise, i.e. we use the universal families $\XX\to B^\circ$, $\YY\to B^\circ$ as in Definition \ref{families}. For each $b\in B^\circ$, there is an isomorphism 
      \begin{equation}\label{isoh0} \Psi_b\colon\ \  t(X_b)\ \xrightarrow{\cong}\ t(Y_b)(-m)\ \ \ \hbox{in}\ \MM_{\rm hom}\ ,\end{equation}
      with an inverse which we will call $\Lambda_b$.

Let  $\pi_X^{\rm var}$
be the projector defining the motives $h_{\rm var}(X)$ as in Proposition \ref{var}. 
On the Pfaffian side, let $h_{\rm var}(Y)$ be the motive defined by the projector
  \[  \pi_Y^{\rm var}:= \Delta_Y -  \sum_{j\not= \dim Y} \pi^j_Y -{1\over d}\,  h^{\dim Y/2}\times h^{\dim Y/2}   \ \ \in\ A^{\dim Y}(Y\times Y)\ ,\]
  where the $\pi^j_Y$ for $j\not=\dim Y$ are the (completely decomposed) projectors given by Proposition \ref{trans}, and $d$ is the degree of $Y$, and it is understood the term $h^{\dim Y/2}$ is zero for $\dim Y$ odd.
  
The nice thing about these projectors is that they are 
{\em generically defined\/} with respect to $B^{\circ}$, i.e. there exist relative cycles
  \[ \pi^{\rm var}_{\XX}\ \ \in\ A^{\dim X_b}(\XX\times_{B^\circ}\XX)\ ,\ \ \  \  \pi^{\rm var}_{\YY}\ \ \in\ A^{\dim Y_b}(\YY\times_{B^\circ}\YY)   \]
  such that
   \[  \pi^{\rm var}_{\XX}\vert_b = \pi^{\rm var}_{X_b}\ ,\ \ \  \ \pi^{\rm var}_{\YY}\vert_b = \pi^{\rm var}_{Y_b}   \ \ \ \ \forall \ b\in B^\circ\ .\]
  Thanks to the condition (NL), we have equalities $t(X)=h_{\rm var}(X)$ and $t(Y)=h_{\rm var}(Y)$ in $\MM_{\rm hom}$ for all fibers over some $B^{\circ\circ}\subset B^\circ$ (where $B^{\circ\circ}$ is a countable intersection of Zariski opens). Combined with the isomorphism \eqref{isoh0}, this means that for all $b\in B^{\circ\circ}$ there are isomorphisms
    \begin{equation}\label{isoh} \Psi_b\colon\ \  h_{\rm var}(X_b)\ \xrightarrow{\cong}\ h_{\rm var}(Y_b)(-m)\ \ \ \hbox{in}\ \MM_{\rm hom}\ ,\end{equation}
      with inverse $\Lambda_b$.

Using Proposition \ref{spread2} (with $\Gamma_\XX=\pi^{\rm var}_\XX$ and $\Gamma_\YY=\pi^{\rm var}_\YY$), we may assume the isomorphism \eqref{isoh} is generically defined, i.e.
there exist relative cycles $\Psi\in A^\ast(\XX\times_{B^\circ}\YY), \Lambda\in A^\ast(\YY\times_{B^\circ}\XX)$ such that 
  \[ \begin{split} ( \Lambda\circ\pi^{\rm var}_\YY\circ \Psi)\vert_b &= \pi^{\rm var}_\XX\vert_b\ \ \ \hbox{in}\ H^{2\dim X_b}(X_b\times X_b,\QQ)\ ,\\
                              ( \Psi\circ\pi^{\rm var}_\XX\circ\Lambda)\vert_b &= \pi^{\rm var}_\YY\vert_b\ \ \ \hbox{in}\ H^{2\dim Y_b}(Y_b\times Y_b,\QQ)\ \ \ \ \ \forall\ b\in B^{\circ\circ}\ .\\
                              \end{split}\]
      
 Applying Proposition \ref{Frtype} (with $M=\bar{M}=\Gr(2,V)$) to the generically defined cycle 
   \[ (\pi^{\rm var}_\XX-\Lambda\circ\pi^{\rm var}_\YY\circ\Psi)\vert_b\ \ \ \in GDA_{hom}^{\dim X_b}(X_b\times X_b)\ ,\] 
 we find that there exists a decomposed cycle $\gamma_b\in A^\ast(X_b)\otimes A^\ast(X_b)$  such that 
   \[ \pi^{\rm var}_{X_b}= ( \Lambda\circ\pi^{\rm var}_\YY\circ\Psi)\vert_b + \gamma_b\ \ \ \hbox{in}\ A^{\dim X_b}(X_b\times X_b)\ .\]
  This means there is a split injection of motives
   \begin{equation}\label{left}  (\Psi_b, \Xi_b)\colon\ \    h_{\rm var}(X_b)\ \hookrightarrow\ h_{\rm var}(Y_b)(-m)\oplus \bigoplus\one(\ast)\ \ \ \hbox{in}\ \MM_{\rm rat}\ \end{equation}
   (with the map from $h_{\rm var}(Y_b)(-m)$ to $h_{\rm var}(X_b)$ given by $\Lambda_b$).
   
   On the Pfaffian side, we can apply Proposition \ref{Frtype} (with $M=\Pf^\circ:=\Pf\setminus \Sing (\Pf)$, which is possible thanks to Example \ref{exPf}; note that the condition $H^\ast_{tr}(Y_b,\QQ) \not=0$ of Proposition \ref{Frtype} is satisfied thanks to Proposition \ref{propcoh}) to the generically defined cycle 
   \[ (\pi^{\rm var}_\YY-\Psi\circ\pi^{\rm var}_\XX\circ\Lambda)\vert_b\ \ \ \in GDA_{hom}^{\dim Y_b}(Y_b\times Y_b)\ .\]    
   The result is that there exists a decomposed cycle $\gamma^\prime_b\in A^\ast(Y_b)\otimes A^\ast(Y_b)$ such that
      \[ \pi^{\rm var}_{Y_b}= ( \Psi\circ\pi^{\rm var}_\XX\circ\Lambda)\vert_b + \gamma^\prime_b\ \ \ \hbox{in}\ A^{\dim Y_b}(Y_b\times Y_b)\ ,\]           
      and hence there is also a split injection of motives
        \begin{equation}\label{right}  (\Lambda_b, \Xi^\prime_b)\colon\ \    h_{\rm var}(Y_b)(-m)\ \hookrightarrow\ h_{\rm var}(X_b)\oplus \bigoplus\one(\ast)\ \ \ \hbox{in}\ \MM_{\rm rat}\ \end{equation}
   (with the map from $h_{\rm var}(X_b)$ to $h_{\rm var}(Y_b)(-m)$ given by $\Psi_b$).

    Taking Chow groups on both sides of \eqref{left} (and exploiting that $t(X)=h_{\rm var}(X)$ has the property that $A^\ast(t(X)) =A^\ast_{hom}(X)$), we find an equality of actions
    \begin{equation}\label{lefta} \Bigl( \pi^{\rm var}_{X_b}\circ \Lambda_b\circ \pi^{\rm var}_{Y_b}\circ \Psi_b\circ \pi^{\rm var}_{X_b} -\pi^{\rm var}_{X_b}  \Bigr){}_\ast=0\colon\ \ \ A^\ast(X_b)\ \to\ A^\ast(X_b)\ \ \ \ \forall\ b\in B^{\circ\circ}\ .\end{equation}
    Likewise, taking Chow groups on both sides of \eqref{right}, we find equality of actions
     \begin{equation}\label{righta} \Bigl( \pi^{\rm var}_{Y_b}\circ \Psi_b\circ \pi^{\rm var}_{X_b}\circ \Lambda_b\circ \pi^{\rm var}_{Y_b} -\pi^{\rm var}_{Y_b}  \Bigr){}_\ast=0\colon\ \ \ 
     A^\ast(Y_b)\ \to\ A^\ast(Y_b)\ \ \ \ \forall\ b\in B^{\circ\circ}\ .\end{equation}
     
   Applying \cite[Proposition 3.2]{SV} (which is inspired by the Bloch--Srinivas decomposition of the diagonal argument \cite{BS}), this means that the expressions in parentheses in \eqref{lefta} and \eqref{righta} are both nilpotent. Taking the largest of the two nilpotence indices, it follows that
     \begin{equation}\label{thisi} \Psi_b\colon\ \ h_{\rm var}(X_b)\ \to\ h_{\rm var}(Y_b)(-m)\ \ \ \hbox{in}\ \MM_{\rm rat} \end{equation}
     is an isomorphism (with inverse $\Lambda^\prime_b$, which is a sum of expressions of the form $\Lambda_b\circ \pi^{\rm var}_{Y_b}\circ \Psi_b\circ\pi^{\rm var}_{X_b}\circ \cdots \circ \Psi_b$), for any $b\in B^{\circ\circ}$.
  In view of the spread lemma \cite[Lemma 3.2]{Vo}, \eqref{thisi} is then an isomorphism for any $b\in B^\circ$.   
 
 This means that for any $b\in B^\circ$, there exist submotives $t(X)\subset h_{\rm var}(X)$, $t(Y)\subset h_{\rm var}(Y)$ as in Proposition \ref{trans}, and such that there is an isomorphism
   \[  t(X_b)\oplus \bigoplus \one(\ast)    \ \xrightarrow{\cong}\ t(Y_b)(-m)\oplus \bigoplus \one(\ast)   \ \ \ \hbox{in}\ \MM_{\rm rat}\ .\]
 Taking Chow groups (and exploiting that $A^\ast(t(X)) =A^\ast_{hom}(X)$), we deduce that there is equality of actions
   \[ \begin{split}   &\Bigl( \pi^{tr}_{X_b}\circ \Lambda^\prime_b\circ \pi^{tr}_{Y_b}\circ \Psi_b\circ \pi^{tr}_{X_b} -\pi^{tr}_{X_b}  \Bigr){}_\ast=0\colon\ \ \ A^\ast(X_b)\ \to\ A^\ast(X_b)   \ ,\\
      &\Bigl( \pi^{tr}_{Y_b}\circ \Psi_b\circ \pi^{tr}_{X_b}\circ \Lambda^\prime_b\circ \pi^{tr}_{Y_b} -\pi^{tr}_{Y_b}  \Bigr){}_\ast=0\colon\ \ \ 
     A^\ast(Y_b)\ \to\ A^\ast(Y_b)\ .  \\   
        \end{split}\]   
 Applying once more \cite[Proposition 3.2]{SV}, this means that the expressions in parentheses are both nilpotent. Taking the largest of the nilpotence indices, we find there is an isomorphism
   \[  \Psi_b\colon\ \ t(X_b)\ \to\ t(Y_b)(-m)\ \ \ \hbox{in}\ \MM_{\rm rat} \ \ \ \ \forall \ b\in B^\circ\ .\]     
   The theorem is proven.
             \end{proof}

 \begin{remark} Concerning the smoothness assumption in Theorem \ref{main}, we add the following observation: if $X$ and $Y$ both have the expected dimension, then $X$ is smooth if and only if $Y$ is smooth. For $(n,k)=(7,7)$ this was observed in \cite{BC}, the general case follows from \cite[Lemma 3.10 and Remark 3.3]{Ji+}.
 \end{remark}

 \section{Applications}
 
 This section presents some applications of the motivic relation of Theorem \ref{main}.
 
 \subsection{Calabi--Yau threefolds} A first consequence concerns the famous Calabi--Yau threefolds of Borisov \cite{Bor}, \cite{BC}:
 
 \begin{corollary}\label{CY3} Let
      \[  \begin{split}   X&:=  \Gr(2,7)\cap   \PP(U) \ \ \ \subset\ \PP^{20} \ , \\
                              Y&:= \Pf\,\cap\, \PP(U^\perp)\ \ \ \subset\ (\PP^{20})^\vee \ 
                              \end{split}  \]  
 be smooth dimensionally transverse intersections, where $\PP(U)\subset\PP^{20}$ is a codimension 7 linear subspace.
 Then $X$ and $Y$ are Calabi--Yau threefolds, and
   \[ h(X)\cong h(Y)\ \ \ \hbox{in}\ \MM_{\rm rat}\ .\]
   \end{corollary}
 
 \begin{proof} Let $h^3(X):= h_{var}(X)$ (which is equal to $t(X)$ because $X$ is odd-dimensional), and similarly for $Y$. According to Theorem \ref{main} with $(n,k)=(7,7)$, there is an isomorphism $h^3(X)\cong h^3(Y)$. Both $X$ and $Y$ have Picard number $1$, and so the result follows.
  \end{proof}

%
%

 \subsection{Cubic threefolds and Fano threefolds of genus 8}
 
 \begin{corollary}\label{65} Let $Y$ be a general cubic threefold. There exists a prime Fano threefold $X$ of genus 8, and an isomorphism of Chow motives
   \[ h(X)\cong h(Y)\ \ \ \hbox{in}\ \MM_{\rm rat}\ .\]
   Moreover, the general prime Fano threefold of genus 8 arises in this way.
   \end{corollary}
   
   \begin{proof} The general cubic threefold is Pfaffian \cite{Be}, and so this is the case $(n,k)=(6,5)$ of Theorem \ref{main}. The ``moreover'' statement follows from work of Mukai \cite{Mu}.
     \end{proof}

   \begin{remark} In the set-up of Corollary \ref{65}, the varieties $X$ and $Y$ are actually birational \cite{Puts}, and the isomorphism of motives can be readily obtained by exploiting the specific form of the birationality \cite{g8}. However, the proof given here does {\em not\/} rely on the birationality.
    \end{remark}

 \subsection{Pfaffian cubic fourfolds}
 
 \begin{corollary}\label{pfafc} Let $Y\subset\PP^5$ be a general Pfaffian cubic fourfold. There exists a genus $8$ K3 surface $X$ and an isomorphism of motives
   \[ h(Y)\cong h(X)(-1)\oplus \one\oplus \one(-4)\ \ \ \hbox{in}\ \MM_{\rm rat}\ .\]
   Moreover, the general genus 8 K3 surface arises in this way.
  \end{corollary}
  
  \begin{proof} This is Theorem \ref{main} with $(n,k)=(6, 6)$. The ``moreover'' statement is work of Mukai \cite{Mu8}. 
   \end{proof}
  
  \begin{remark} Pfaffian cubic fourfolds and genus $8$ K3 surfaces are also related on the level of derived categories \cite[Theorem 10.4]{Kuz0}. The relation of motives of Corollary \ref{pfafc} is well-known, and is also valid for other cubic fourfolds with an associated K3 surface \cite[Theorem 0.3]{Bul}. Contrary to \cite{Bul}, however, the argument proving Corollary \ref{pfafc} is direct and geometric, and does not rely on derived category results.
  \end{remark}

\subsection{Varieties with trivial Chow groups}

The following is not a corollary of the main result. However, we include it because it fits in well here:

\begin{proposition}\label{12} Let $X$ be a smooth dimensionally transverse intersection
   \[ X:= \Gr(2,n)\cap H_1\cap\cdots\cap H_k\ \ \ \subset\ \PP^{{n\choose 2}-1} \ ,\]
   where $H_j$ are hyperplanes. Assume $k\le 2$. Then
     \[ A^i_{hom}(X)=0\ \ \ \forall\ i\ .\]
\end{proposition}

\begin{proof} It is well-known that $H^\ast_{tr}(X,\QQ)=0$ in this case \cite[Propositions 2.3 and 2.4]{Don}, and so Theorem \ref{main} does not apply. However, one can understand the Chow groups of $X$ by using a straightforward geometric argument, inspired by \cite{Don}.

Let $P\subset\PP(V_n)$ be a fixed hyperplane, and consider (as in \cite[Section 2.3]{Don}) the rational map
  \[  \Gr(2,V_n)\ \dashrightarrow\ P \]
  sending a line in $\PP(V_n)$ to its intersection with $P$. This map is resolved by blowing up a subvariety $\sigma_{11}(P)\cong \Gr(2,n-1)$, resulting in a morphism
   \[ \Gamma\colon\ \ \wt{\Gr}\ \to\ P \]
   (where $\wt{\Gr}\to \Gr(2,V_n)$ denotes the blow-up with center $\sigma_{11}(P)$).

Let $\wt{X}\to X$ be the blow-up of $X$ with center $\sigma_{11}(P)\cap X$, and let us consider the morphism
  \[ \Gamma_X\colon\ \ \wt{X}\ \to\ P \ ,\]
  obtained by restricting $\Gamma$.
  
  In case $k=1$ and $P$ is generic with respect to $X$, the morphism $\Gamma_X$ is a $\PP^{n-3}$-fibration over $P$. It follows that $\wt{X}$, and hence $X$, has trivial Chow groups.
  
  In case $k=2$, and $P$ chosen generically with respect to $X$, the morphism $\Gamma_X$ is generically a $\PP^{n-4}$-fibration over $P$, and there are finitely many points in $P$ where the fiber is $\PP^{n-3}$. Applying Theorem \ref{ji}, this implies that $\wt{X}$, and hence $X$, has trivial Chow groups.
\end{proof}     

\begin{remark} In case $n$ is 6 or 7, and $X\subset\Gr(2,n)$ as in Proposition \ref{12}, it follows from HPD that the derived category of $X$ has a full exceptional collection
\cite[Corollaries 10.2 and 11.2]{Kuz0}. Since it is known that varieties admitting a full exceptional collection have trivial Chow groups \cite{MT}, this gives another proof of Proposition \ref{12} for $n=6$ or $n=7$.
\end{remark}

 \subsection{Fano varieties with finite-dimensional motive}
 
 \begin{corollary}\label{findim} Let $X$ be a smooth dimensionally transverse intersection
   \[ X:= \Gr(2,n)\cap H_1\cap\cdots\cap H_k\ \ \ \subset\ \PP^{{n\choose 2}-1} \ ,\]
   where $H_j$ are hyperplanes.
   Assume either $k=3$ and $n$ is even, or $k=5$ and $n$ is odd, or $(n,k)=(8,5)$.
   Then
     \[ A^\ast_{AJ}(X)=0\ ,\]
     and in particular
      $X$ has finite-dimensional motive (in the sense of \cite{Kim}).
     \end{corollary}
     
     \begin{proof} We first treat the cases $k=3$ and $n$ is even, or $k=5$ and $n$ is odd.
     As a first step, let us assume $X$ is such that the dual $Y\subset (\PP^{{n\choose 2}-1})^\vee $ is smooth and dimensionally transverse, i.e. $Y$ has dimension $1$.
     Theorem \ref{main} implies the isomorphism
      \[  h(X)\cong   t(Y)(-m)\oplus \bigoplus\one(\ast)\ \ \ \hbox{in}\ \MM_{\rm rat}\ .\]
     Since curves have finite-dimensional motive and injective Abel--Jacobi maps, this implies both statements for $X$.
     
     Next, let $\XX\to B$ denote the universal family of all smooth complete intersections of the type under consideration, and let $B^\circ\subset B$ denote the Zariski open subset parametrizing smooth $X$ for which the dual $Y$ is a smooth curve. 
%
  By definition of finite-dimensionality, the above means exactly that for all $b\in B^\circ$ one has vanishing
    \begin{equation}\label{find}   \sym^{b_{tr}} t(X_b) =0\ \ \      \hbox{in}\ A^{b_{tr}\cdot\dim X_b}(X_b^{b_{tr}}\times X_b^{b_{tr}})\  ,    \end{equation}
   where the symmetric product of a motive is as in \cite{Kim}, and
  $b_{tr}:=\dim H^{\dim X_b}(X_b,\QQ)$.
  But the projector defining $t(X_b)$ is generically defined, and so the spread lemma \cite[Lemma 3.2]{Vo} then implies that the vanishing \eqref{find} holds for all $b\in B$, i.e. all $X_b$ have finite-dimensional motive.      
  
  To prove that $A^\ast_{AJ}(X_b)=0$ for all $X_b$, we observe that the above implies that for all $b\in B^\circ$ one has
   \begin{equation}\label{niv} \hbox{Niveau}\bigl(A^\ast(X_b)\bigr)\le 1    \end{equation}
   in the sense of \cite{moi}. This means that for each $b\in B^\circ$
    there is a decomposition of the diagonal
     \[ \Delta_{X_b} = \gamma_b^0+\cdots +\gamma_b^r\ \ \ \hbox{in}\ A^{\dim X_b}(X_b\times X_b)\ \]
     where $\gamma_b^j$ is supported on $V_b^j\times W_b^j\subset X_b\times X_b$ and $\dim V_b^j+\dim W_b^j=\dim X_b+1$.    
     Using the Hilbert schemes argument of \cite[Proposition 3.7]{V0} (cf. also \cite[Proposition A.1]{LNP} for the precise form used here), the 
     $\gamma_b^j, V_b^j, W_b^j$ exist relatively, i.e.
     one can find subvarieties $\VV^j, \WW^j\subset \XX$ with $\codim \VV_j+\codim \WW_j =\dim X_b-1$, and cycles $\gamma^j$ supported on $\VV^j\times_{B^\circ} \WW^j$ such that
     \[ \Delta_\XX\vert_b= \gamma^0\vert_b +\cdots + \gamma^r\vert_b\ \ \ \hbox{in}\ A^{\dim X_b}(X_b\times X_b)\ \ \ \forall\ b\in B^\circ\ .\ \]
     Let $\bar{\gamma}^j\in A^{\dim X_b}(\XX\times_B \XX)$ be cycles that restrict to $\gamma^j\in A^{\dim X_b}(\XX\times_{B^\circ} \XX)$, and let $\bar{\VV}^j\times_B \bar{\WW}^j$ be the support of $\bar{\gamma}^j$. Given $b_1\in B\setminus B^\circ$, it may happen that $\dim \bar{\VV}^j\vert_{b_1} +  \dim  \bar{\WW}^j\vert_{b_1} $ is larger than $\dim X_{b_1}+1$. However, using the moving lemma, one can find representatives for $\bar{\gamma}^j$ such that the supports verify the condition 
     \[ \dim \bar{\VV}^j\vert_{b_1} +  \dim  \bar{\WW}^j\vert_{b_1} = \dim X_{b_1}+1\ ,\]
     i.e. \eqref{niv} holds for $X_{b_1}$. 

We have ascertained that \eqref{niv} holds for all $b\in B$.
Letting the decomposition of the diagonal act on Chow groups, this shows that
     \[  A^\ast_{AJ}(X_b)=0\ \ \ \forall\ b\in B\ .\]

 Finally, the argument for $(n,k)=(8,5)$ is similar: in this case, $X$ is a Fano sevenfold and $Y$ is a Fano threefold. Since Fano threefolds have finite-dimensional motive, Theorem \ref{main} implies the same for $X$.    
      \end{proof}

  \begin{remark} Corollary \ref{findim} improves on results of Tabuada, who established (using HPD results and non-commutative motives) Schur-finiteness for linear sections of $\Gr(2,5)$,
  $\Gr(2,6)$ and $\Gr(2,7)$ \cite[Theorems 1.3 and 1.4]{Tab}.
  \end{remark}

%
%
%

 \subsection{Fano varieties satisfying Voevodsky's conjecture}
 
 \begin{definition}[Voevodsky \cite{Voe}]\label{sm} Let $X$ be a smooth
		projective variety. A cycle $a\in A^i(X)$ is called {\em smash-nilpotent\/} 
		if there exists $m\in\NN$ such that
		\[ \begin{array}[c]{ccc}  a^m:= &\undermat{(m\hbox{
				times})}{a\times\cdots\times a}&=0\ \ \hbox{in}\ 
		A^{mi}(X\times\cdots\times
		X)_{}\ .
		\end{array}\]
		\vskip0.6cm
		
		Two cycles $a,a^\prime$ are called {\em smash-equivalent\/} if their
		difference
		$a-a^\prime$ is smash-nilpotent. We will write $A^i_\otimes(X)\subseteq
		A^i(X)$ for the subgroup of smash-nilpotent cycles.
	\end{definition}
	
	\begin{conjecture}[Voevodsky \cite{Voe}]\label{voe} Let $X$ be a smooth
		projective variety. Then
		\[  A^i_{num}(X)\ \subseteq\ A^i_\otimes(X)\ \ \ \hbox{for\ all\ }i\ .\]
	\end{conjecture}
	
	\begin{remark} It is known \cite[Th\'eor\`eme 3.33]{An} that Conjecture
		\ref{voe} for all smooth projective varieties implies (and is strictly
		stronger than) Kimura's conjecture ``all smooth projective varieties have finite-dimensional motive'' \cite{Kim}.
	\end{remark}

  \begin{corollary}\label{smash} Let $X$ be a smooth dimensionally transverse intersection
   \[ X:= \Gr(2,n)\cap H_1\cap\cdots\cap H_k\ \ \ \subset\ \PP^{{n\choose 2}-1} \ ,\]
   where $H_j$ are hyperplanes and $k\in\{3,4\}$ if $n$ is even and $k\in\{5,6\}$ if $n$ is odd.
   Then Voevodsky's conjecture is true for $X$. Moreover,
     \[ A^j_{AJ}(X)=0\ \ \ \forall j\not= {1\over 2}\dim X+1\ .\]
     \end{corollary}
  
  \begin{proof} In a first step, we assume $X$ is such that the dual $Y\subset\Pf$ on the Pfaffian side is smooth and dimensionally transverse. Theorem \ref{main} applies
   and gives a split injection 
     \[ h(X)\ \hookrightarrow\ h(Y)(-m)\oplus \bigoplus\one(\ast)\ \ \ \hbox{in}\ \MM_{\rm rat}\ .\]
     Since $\dim Y\le 2$, this implies that
     \[ \hbox{Niveau}\bigl(A^\ast(X)\bigr)\le 2 \]
     (in the sense of \cite{moi}), i.e. there is a decomposition of the diagonal
     \begin{equation}\label{dd} \Delta_X = \gamma_0+\cdots +\gamma_r\ \ \ \hbox{in}\ A^{\dim X}(X\times X)\ \end{equation}
     where $\gamma_j$ is supported on $V_j\times W_j\subset X\times X$ and $\dim V_j+\dim W_j=\dim X+2$.  
     Looking at the action on Chow groups, this decomposition implies that homological and algebraic equivalence coincide on $X$. Since it is known that $A^\ast_{alg}(X)\subset A^\ast_{\otimes}(X)$ \cite{Voe}, \cite{V9}, this implies that Conjecture \ref{voe} holds for $X$.
     
 Next, let us extend to arbitrary $X$. Let $\XX\to B$ denote the universal family, and $B^\circ\subset B$ the Zariski open where $X$ and its dual $Y$ are simultaneously smooth and dimensionally transverse. The first step above shows that $ \hbox{Niveau}\bigl(A^\ast(X_b)\bigr)\le 2 $ for all $b\in B^\circ$. As in the proof of Corollary \ref{findim}, this property extends from $B^\circ$ to $B$ and we find that 
     \[ \hbox{Niveau}\bigl(A^\ast(X_b)\bigr)\le 2 \ \ \ \forall\ b\in B\ .\]   
     In particular, Voevodsky's conjecture holds for all $X_b$ in the family.
     
     Finally, the vanishing of $A^j_{AJ}(X)$, $j\not= {1\over 2}\dim X+1$ is a straightforward consequence of the decomposition \eqref{dd}.     
      \end{proof}   
     
 \begin{remark} In \cite[Theorem 1.7]{BMT}, the special case of Corollary \ref{smash} where $n=6$ or $n=7$ was proven (in loc. cit., the restriction to $n=6,7$ is necessary because the argument relies on HPD for Grassmannians via non-commutative motives). The argument proving Corollary \ref{smash} is more elementary, in that we do not rely on derived category arguments at all.
 \end{remark}

 \subsection{Fano varieties satisfying the Hodge conjecture}
 
 \begin{corollary}\label{hc} Let $X$ be a smooth dimensionally transverse intersection
   \[ X:= \Gr(2,n)\cap H_1\cap\cdots\cap H_k\ \ \ \subset\ \PP^{{n\choose 2}-1} \ ,\]
   where $H_j$ are hyperplanes and either $n$ is even and $k\in\{3,4,5\}$ or $n$ is odd and $k\in\{5,6\}$. Then
      \[ A^j_{hom}(X)=0\ \ \ \forall \ j> {1\over 2}(\dim X+3)\ ,\]
  and in particular the Hodge conjecture is true for $X$.
   \end{corollary}  
   
   \begin{proof} First, let us assume $X$ is such that the dual $Y\subset\Pf$ on the Pfaffian side is smooth and dimensionally transverse. Theorem \ref{main} applies
   and gives a split injection 
     \[ h(X)\ \hookrightarrow\ h(Y)(-m)\oplus \bigoplus\one(\ast)\ \ \ \hbox{in}\ \MM_{\rm rat}\ .\]
     Since $\dim Y\le 3$, this implies that
     \[ \hbox{Niveau}\bigl(A^\ast(X)\bigr)\le 3 \]
     (in the sense of \cite{moi}), i.e. there is a decomposition of the diagonal
     \[ \Delta_X = \gamma_0+\cdots +\gamma_r\ \ \ \hbox{in}\ A^{\dim X}(X\times X)\ \]
     where $\gamma_j$ is supported on $V_j\times W_j\subset X\times X$ and $\dim V_j+\dim W_j=\dim X+3$.
     
     Next, let us extend to arbitrary $X$. Let $\XX\to B$ denote the universal family, and $B^\circ\subset B$ the Zariski open where $X$ and its dual $Y$ are simultaneously smooth and dimensionally transverse. The above shows that $ \hbox{Niveau}\bigl(A^\ast(X_b)\bigr)\le 3   $ for all $b\in B^\circ$. As in the proof of Corollary \ref{findim}, this property extends from $B^\circ$ to $B$ and we find that 
     \[ \hbox{Niveau}\bigl(A^\ast(X_b)\bigr)\le 3 \ \ \ \forall\ b\in B\ .\]   
     It is well-known (cf. for instance \cite{moi}) that this implies the vanishing of $A^j_{hom}(X)$ in the indicated range, as well as the truth of the Hodge conjecture for $X_b$.
         \end{proof}
         
    In the special case $(n,k)=(6,3)$ one can say more:
         
     \begin{corollary}\label{hc2} Let $X$ be a general intersection
   \[ X:= \Gr(2,6)\cap H_1\cap H_2 \cap H_3\ \ \ \subset\ \PP^{14} \ ,\]
   where $H_j$ are hyperplanes. Then $X$ is a Fano fivefold, and the generalized Hodge conjecture is true for $X^m$ for all $m\in\NN$.
   \end{corollary}
   
   \begin{proof} $X$ being general, the dual $Y\subset\Pf$ on the Pfaffian side is an elliptic curve. The isomorphism of motives of Theorem \ref{main} implies there is an isomorphism of Hodge structures
    \[ H^j(X^m,\QQ)\ \cong\  H^{j-4m}(Y^m,\QQ)(-2m) \oplus \bigoplus H^\ast(Y^{m-1},\QQ)(\ast)\oplus  \cdots \oplus \bigoplus \QQ(\ast)\ .\]
   Since this isomorphism is also compatible with the coniveau filtration \cite[Proposition 1.2]{V4}, one is reduced to proving the generalized Hodge conjecture for powers of an elliptic curve $Y$. This is known thanks to work of Abdulali \cite[Section 8.1]{Ab} (cf. also \cite[Corollary 3.13]{Vial}). 
   
   \end{proof}         
 
 \subsection{Fano eightfolds of K3 type} 
 
  \begin{corollary}\label{8f} Let $X$ be a general complete intersection
   \[ X:= \Gr(2,8)\cap H_1\cap H_2\cap H_3\cap H_4\ \ \ \subset\ \PP^{ 27 } \ ,\]
   where $H_j$ are hyperplanes. Then there exists a quartic K3 surface $S$ and an isomorphism of motives
     \[ h(X)\cong h(S)(-3)\oplus \bigoplus\one(\ast)\ \ \ \hbox{in}\ \MM_{\rm rat}\ .\]
      \end{corollary} 
      
     \begin{proof} This is the case $(n,k)=(8,4)$ of Theorem \ref{main}.
      \end{proof}
     
 
 \begin{corollary}\label{rho} For any $3\le \rho\le 22$, there exist Fano eightfolds $X$ as in Corollary \ref{8f} with
 \[   \dim \ima\Bigl( A^4(X)\to H^8(X,\QQ)\Bigr) = \rho\ .\]
 For $\rho\ge 21$, $X$ has finite-dimensional motive.
 \end{corollary}
 
 \begin{proof} 
 Inside the moduli space $\FF_3$ of genus 3 K3 surfaces, let $\FF_3^\circ \subset \FF_3$ denote the locus of K3 surfaces which are Pfaffian quartics and for which one of the duals $X\subset\Gr(2,8)$ is smooth and dimensionally transverse.
 As the general quartic K3 surface is Pfaffian \cite{Be}, $\FF_3^\circ$ is a dense open subset. Let $\FF_3^{\ge \rho}\subset\FF_3$ denote the locus parametrizing K3 surfaces with Picard number $\ge\rho$.
  The Noether--Lefschetz theory for K3 surfaces \cite[Chapter 17]{Huyb} implies that the locus
 $\FF_3^{\ge \rho}$ has dimension $20-\rho$ and is analytically dense in $\FF_3^{\ge \rho-1}$, and so in particular for each $1\le \rho\le 20$, the locus $\FF_3^{\ge \rho}$ meets the open
 $\FF_3^\circ$.
 Given $Y$ a K3 surface in $\FF_3^{\ge \rho}\cap \FF_3^\circ$, let $X\subset\Gr(2,8)$ be a dual eightfold.
  The isomorphism of Corollary \ref{8f} induces an isomorphism
   \[  H^8_{var}(X,\QQ)\cap \ima\Bigl( A^4(X)\to H^8(X,\QQ)\Bigr)  \cong H^2_{var}(Y,\QQ)\cap  \ima\Bigl( A^1(Y)\to H^2(Y,\QQ)\Bigr) \ ,\]
   and so 
         \[   \dim \ima\Bigl( A^4(X)\to H^8(X,\QQ)\Bigr) = \rho+2\ .\]   
         
  The last statement follows from Corollary \ref{8f} plus the fact that K3 surfaces with Picard number $\ge 19$ have finite-dimensional motive \cite{Ped}.       
          \end{proof}
     
  \begin{remark} In particular, Corollary \ref{rho} with $\rho=22$ gives examples of $\rho$-maximal varieties, i.e. $2m$-dimensional varieties $X$ with the property that
    \[   \dim \ima\Bigl( A^m(X)\to H^{2m}(X,\QQ)\Bigr) = \dim H^{m,m}(X,\C)  \ .\]
    This notion is discussed in \cite[Section 8]{Beau}.
    \end{remark}

 \subsection{Fano varieties of Calabi--Yau type with infinite-dimensional Griffiths group} 
 
 \begin{definition} Let $X$ be a smooth projective variety. The {\em Griffiths groups\/} of $X$ are defined as
   \[ \Grif^j(X):= {z^j_{hom}(X) \over z^j_{alg}(X)}\ ,\]
   where $z^j_{hom}(X)$ and $z^j_{alg}(X)$ denote the groups of codimension $j$ algebraic cycles on $X$ that are homologically trivial, resp. algebraically trivial.
   \end{definition}

 \begin{corollary}\label{grif} Let $X$ be a general complete intersection
   \[ X:= \Gr(2,10)\cap H_1\cap\cdots\cap H_5\ \ \ \subset\ \PP^{44} \ ,\]
   where the $H_j$ are Pl\"ucker hyperplanes.
 Then $X$ is a Fano elevenfold and the Griffiths group $\Grif^6(X)_{\QQ}$ is infinite-dimensional.  
  \end{corollary}
  
  \begin{proof} For sufficiently general hyperplanes $H_j$, both $X$ and its dual 
     \[ Y:= \Pf \,\cap\, \PP(U^\perp)\ \ \ \subset\ (\PP^{44})^\vee \ \]
     are smooth and dimensionally transverse. Theorem \ref{main} with $(n,k)=(10,5)$ then gives an isomorphism of motives
     \[ h(X)\cong t(Y)(-4)\oplus \bigoplus\one(\ast)\ \ \ \hbox{in}\ \MM_{\rm rat}\ .\]
     Taking Griffiths groups, this implies
     \[  \Grif^j(X)_{\QQ}=0\ \ \forall\ j\not=6\ ,\ \ \ \Grif^6(X)_{\QQ}\cong \Grif^2(Y)_{\QQ}\ .\]
     Here $Y$ is a quintic threefold, and the general quintic threefold arises in this way \cite[Proposition 8.9]{Be}. The corollary thus follows
     from Clemens' celebrated result that $\Grif^2(Y)_{\QQ}$ is infinite-dimensional for a general quintic threefold $Y$ \cite{Clem}.
     \end{proof}

 \begin{corollary}\label{grif2} Let $Z$ be an intersection
  \[ Z:= \Bigl( \Gr(2,10)\times \PP^4\Bigr)\cap H_{(1,1)}\ \ \ \subset \PP^{44}\times\PP^4\ ,\]
  where $H_{(1,1)}$ is a general bidegree $(1,1)$ hypersurface. Then $Z$ is a Fano 19-fold and the Griffiths group $\Grif^{10}(Z)_{\QQ}$ is infinite-dimensional.   
 \end{corollary}
 
 \begin{proof} Using the Cayley trick in the form of Theorem \ref{ji}, we find there is an isomorphism of motives
   \[ h(Z)\oplus \bigoplus \one(\ast)\cong h(X)(-4)\oplus \bigoplus \one(\ast)\ \ \ \hbox{in}\ \MM_{\rm rat}\ ,\]
   where $X$ is an intersection of $\Gr(2,10)$ with 5 general hyperplanes. Taking Griffiths groups, this implies in particular
    \[  \Grif^j(Z)_{\QQ}=0\ \ \forall\ j\not=10\ ,\ \ \ \Grif^{10}(Z)_{\QQ}\cong \Grif^6(X)_{\QQ}\ .\]
    The corollary now follows from Corollary \ref{grif}.
 \end{proof}
 
 \begin{remark} Following up on Clemens' famous result about the quintic threefold, infinite-dimensionality of the Griffiths group has been proven in \cite{AC} for the cubic sevenfold, and in \cite{FIK} for certain other Fano varieties of Calabi--Yau type. The results in \cite{AC} and \cite{FIK}
 are more difficult and remarkable than those of Corollary \ref{grif} and \ref{grif2}: indeed, the varieties of \cite{AC} and \cite{FIK} correspond to a non-commutative Calabi--Yau threefold (i.e. the interesting part of the derived category is a CY3 category without geometric incarnation), while the varieties of Corollary \ref{grif} and \ref{grif2} are (motivically and categorically) related to an honest Calabi--Yau threefold.
 \end{remark}

 \vskip1cm
\begin{nonumberingt} Thanks to Lie Fu, Mingmin Shen and Charles Vial for sundry inspiring exchanges in Amsterdam in February 2020. Thanks to the referee for many helpful and pertinent comments that significantly improved the paper. Thanks to Kai and Len for being assiduous Tintin readers.
\end{nonumberingt}

\end{document}